\documentclass[final]{siamltex}

\usepackage{amssymb,amsmath,epsfig}

\title{The Relation between  Approximation in Distribution and Shadowing  in Molecular Dynamics\thanks{This work was supported by a National Science and Engineering Research Council Discovery Grant.}}

\author{Paul Tupper\thanks{Department of Mathematics, Simon Fraser University, 8888 University Drive, Burnaby BC, V5A 1S6, Canada.({\tt pft3@math.sfu.ca }).}}

\begin{document}

\maketitle

\begin{abstract}
Molecular dynamics refers to the computer simulation of a material at the atomic level.
An open problem in numerical analysis is to explain the apparent reliability of molecular dynamics simulations.
 The difficulty is that individual trajectories computed in molecular dynamics are accurate for only short time intervals, whereas
apparently reliable information can be extracted from very long-time simulations.
  It has been conjectured that
 long molecular dynamics trajectories  have low-dimensional statistical features that accurately approximate those of the original system.
Another conjecture is that numerical trajectories satisfy the shadowing property: that they are close over long time intervals to exact trajectories but with different initial conditions. We prove that these two views are actually equivalent to each other, after we suitably modify the concept of shadowing.   A key ingredient of our result is a general theorem that  allows us to take random elements of a metric space that are close in distribution and embed them in the same probability space so that they are close in a strong sense.
This result is similar to the Strassen-Dudley Theorem except that a mapping is provided between the two random elements.
Our results on shadowing are motivated by molecular dynamics but apply to the approximation of any dynamical system when initial conditions are selected according to a probability measure.
\end{abstract}


\pagestyle{myheadings}
\thispagestyle{plain}
\markboth{P. F. TUPPER}{APPROXIMATION IN DISTRIBUTION AND SHADOWING}

\section{Introduction} \label{sec:intro}

In the form we consider here, molecular dynamics
consists of modeling an atomistic system with a system of Hamiltonian ordinary differential equations that are numerically integrated.   For given initial conditions information of physical interest is extracted from the resulting approximate trajectories \cite{frenkel,allen,LeiRei}.
Despite the scientific importance of molecular dynamics there is very little rigorous justification of the results it produces.
The problem
is that individual trajectories computed by molecular dynamics
simulations are accurate for only short time intervals.   Numerical trajectories diverge
rapidly from true trajectories given the step-lengths used in
practice.  
This phenomenon is well-known but is not considered a short-coming of molecular dynamics in practical terms \cite[p.\ 81]{frenkel}.
Experience has suggested that the features of trajectories that researchers wish to study are computed accurately.  
  However, that trajectories are reliable in this sense has yet to be rigorously
demonstrated in representative cases. 

Two proposals have emerged to explain the success of molecular dynamics given the inaccuracy of the computed trajectories \cite{SkeTup}.   The first we refer to as \emph{approximation in distribution} and the second as \emph{shadowing}.

The idea of approximation in distribution is to view both exact trajectories and numerical trajectories as stochastic processes.  This is done by drawing initial conditions from a physically appropriate probability distribution rather than considering a single fixed initial condition.
 Both the resulting numerical trajectory and the resulting exact trajectory are then random.  The proposal is that in some distributional (statistical) sense the numerical trajectories and exact trajectories are close to each other.
This means roughly that the probability of some event happening for the numerical trajectory is close to that of the same event happening for the exact trajectory. 
Put in this way, approximation in distribution may not hold for arbitrary events
which depend on the position of every single atom in the system.
 However, usually we are only interested in low-dimensional functions of the state of the full system.     An example we shall consider in the following is when one is only interested in the position of a single particle in a system consisting of many particles.
It may be that for systems studied in molecular dynamics statistical features of trajectories of single particles are reproduced accurately in simulations.
Approximation in distribution for low dimensional functions of numerical trajectories has been studied for model systems in \cite{cano,hald,tupperima} using a combination of analysis and computational experiments, though it has not been established for more realistic systems.

 The idea of shadowing is to show that, even though a numerical trajectory diverges rapidly from the corresponding exact trajectory, it may be possible to show that the numerical trajectory is close to another exact trajectory with different initial conditions.
 If shadowing holds, then numerical trajectories from molecular dynamics simulations can be viewed as real trajectories with some small observational error.
Shadowing has been established for various types of dynamical systems with uniform hyperbolicity properties \cite{pilyugin} and these ideas have been applied to Hamiltonian dynamical system such as those studied in molecular dynamics \cite{Reich}. 
However, shadowing over arbitrarily long time intervals is probably not possible for 
realistic Hamiltonian systems \cite{saueryorke}.
Moreover, shadowing as it is usually defined does not guarantee that statistical features of numerical trajectories match those of the exact trajectories \cite{Hayes}.
We will discuss these issues further and 
show how to modify the concept of shadowing suitably so that it is applicable to our case and to other situations where the initial conditions of the dynamical system are distributed according to some probability measure.

The main purpose of this paper is to carefully define and quantify these two concepts, approximation in distribution and shadowing, in the context of molecular dynamics and to explain the relation between them.
 Our main result shows that when the two concepts  are formalized and suitably modified they are actually equivalent.

In Section~\ref{sec:problem} we introduce a model system for molecular dynamics and present the results of some numerical experiments performed with it.
We first demonstrate that numerical trajectories using practical time steps diverge rapidly from exact trajectories.  We then provide evidence that statistical features of some low-dimensional functions of the trajectories are nevertheless reliable.

In Section~\ref{sec:approaches} we discuss approximation in distribution and shadowing in detail and give quantitative versions of each idea.  
 In particular we show how to adapt the idea of shadowing to situations where initial conditions are distributed according to a probability measure.
 Our concept, which is a modification of the usual notion of shadowing for dynamical systems,  we call Weak Shadowing.

In Section~\ref{sec:mainthm} we prove our main result, Theorem~\ref{thm:easythm}, which we state here.
  In what follows, the space $(C[0,T])^m$ is the set of all continuous trajectories on $[0,T]$ taking values in $\mathbb{R}^m$.  For $x,y \in (C[0,T])^m$, we define $\|x -y \|_\infty = \sup_{t \in [0,T]} | x(t) -y(t) |$.   We let $\Pi$ be a map $\mathbb{R}^m \rightarrow \mathbb{R}^k$ and we define $\Pi(x) \in (C[0,T])^k$ by $\Pi(x)(t) = \Pi(x(t))$ for any $x \in (C[0,T])^m$.
 We let $X_0$ be a random initial condition in $\mathbb{R}^m$ and   then denote by $X$ the   random member of $(C[0,T])^m$ starting at $X_0$ and generated by the differential equations.
  When we use a numerical method to generate an approximate solution to the differential equations at a sequence of points, we use $X_{\Delta T}$ to denote the random element of $(C[0,T])^m$ generated by its linear interpolation.
Finally,  we denote by $\rho$ the well-known Prokhorov metric on random elements of metric spaces which we will define in Subsection~\ref{subsec:distrib}.   It metrizes convergence in distribution so that if two random elements of a metric space are close according to $\rho$ they have approximately the same distribution. 

\begin{theorem} \label{thm:easythm}
Let $X_0$ be a random vector in $\mathbb{R}^m$ such that $\mathbb{P}(X_0 =x)=0$ for any $x \in \mathbb{R}^m$. Let $X$ be the random trajectory in $(C[0,T])^m$ generated by a system of differential equations starting from $X_0$.   Let $X_{\Delta t}$ be the random trajectory of $(C[0,T])^m$ generated by a numerical method starting at $X_0$. 
Let $\Pi : \mathbb{R}^m \rightarrow \mathbb{R}^k$ be a map.
  Then the following are equivalent for all $\epsilon>0$: \\
(A) {\bf Approximation in distribution}.
\[
\rho( \Pi(X), \Pi(X_{\Delta t}) )  < \epsilon.
\]
(B) {\bf Weak Shadowing}.  There is a map $\mathcal{S}_{\Delta t} \colon \mathbb{R}^m \rightarrow \mathbb{R}^m$ such that $Y_0 = \mathcal{S}_{\Delta t} X_0$ has the same probability distribution as $X_0$ and 
if $Y$ is the random member of $(C[0,T])^m$ starting at $Y_0$ generated by the flow of the differential equations then
\[
\mathbb{P}(  \|  \Pi(X_{\Delta t}) - \Pi(Y)  \|_{\infty} > \epsilon ) < \epsilon.
\]
\end{theorem}

Note that no assumptions are made about the differential equations that generate the exact trajectory $X$ nor about the numerical method that generates the approximate trajectory $X_{\Delta}$.  The theorem just asserts that two ways in which a numerical method can be accurate, (A) and (B), are equivalent.  The theorem does not assert that (A) or (B) holds for any particular system or any particular method.

   Showing that (B) implies (A) is straightforward, but the converse requires the result of Theorem~\ref{thm:maintechnical}, which is a version of the Strassen-Dudley theorem \cite{strassen}, \cite[$\S$11.6]{dudley}.  
 The original Strassen-Dudley Theorem shows that two random variables that are close with respect to the Prokhorov metric $\rho$ can be embedded in a new probability space where they are close in a strong sense.   
  Our  contribution is to show how to do this with one random variable defined as a function of the other.
The only important extra assumption needed is that the measures induced by the random variables be non-atomic, which means that they assign zero measure to any point.
     
  Finally, in Section~\ref{sec:discussion} we conclude with a discussion of what our result suggests for the numerical analysis of molecular dynamics.

\section{Numerical Experiments} \label{sec:problem}

We consider a  system of $n=100$ point particles interacting on an 11.5
by 11.5 square periodic domain.  We let $q \in \mathbb{T}^{2n}$ and $p
\in \mathbb{R}^{2n}$ denote the positions and velocities of the
particles, with $q_i \in \mathbb{T}^2, p_i \in \mathbb{R}^2$ denoting
the position and velocity of particle $i$.   The motion of the particles
is described by a system of Hamiltonian differential equations: 
\begin{equation} \label{eqn:hamdynamics}
\frac{dq}{dt} = \frac{\partial H}{\partial p}, \ \ \ 
\frac{dp}{dt} = -\frac{\partial H}{\partial q},
\end{equation}
with Hamiltonian 
\[
H(q,p) = \frac{1}{2} \|p\|_2^2 + \sum_{i<j} V_{LJ} ( \| q_i-q_j \| ).
\]
Here $V_{LJ}$ denotes the famous Lennard-Jones potential \cite{frenkel}.  In our
simulations we use a truncated but infinitely smooth version \cite[p.\ 2409]{StillWeb}: 
\[
V_{LJ}(r) = \left\{  
\begin{array} {ll}
4 \left( \frac{1}{r^{12}} - \frac{1}{r^6} \right) \exp[(r-r_{\mbox{\tiny{cutoff}}})^{-1}] , & \mbox{if } r
\leq r_{\mbox{\tiny{cutoff}}}, 
\\
0,& \mbox{otherwise.}
\end{array} \right.
\]
We set $r_{\mbox{\tiny{cutoff}}} = 2.5$.

For our first numerical experiments we take our initial conditions $q(0),p(0)$ to be randomly distributed
 according to the probability density function 
\begin{equation} \label{eq:distrib} 
 C e^{-\beta H(q,p)},
\end{equation}
where $C$ is chosen so that $\int C e^{-\beta H(q,p)} \, dq \, dp = 1$.   
The probability distribution with this density function is known as the canonical distribution for the system with Hamiltonian $H$ at temperature $\beta^{-1}$.    It is intended to model the equilibrium distribution of the system when it is in thermal contact with an environment of temperature $\beta^{-1}$ \cite[Sec.\ 6.2]{reif}.
    For our experiments we fixed $\beta=1$ and generate $(q(0), p(0))$ according to \eqref{eq:distrib} using Langevin dynamics \cite{cances}.  We then subtract a constant vector from the velocities of all the particles so that center of mass of the system has zero velocity.
The canonical distribution for any $\beta >0$ with this adjustment  is invariant with respect to the dynamics described by \eqref{eqn:hamdynamics}.  Later in this section we will perform further experiments with a nonequilibrium distribution on the initial conditions.
 
We numerically integrate \eqref{eqn:hamdynamics} using the
St\"ormer-Verlet scheme, which is an explicit second-order method for our system \cite{hairer}.  It is the standard numerical integrator used in molecular dynamics \cite[p.\ 69]{frenkel}.
Given an initial $(q^0, p^0)=(q(0),p(0))$ and a $\Delta t>0$, the St\"ormer-Verlet scheme generates a sequence
of states $(q^n, p^n)$, $n\geq0$ such that $(q^n,p^n) \approx (q(n\Delta
t),p(n\Delta t))$.  The version of the algorithm we use is  
\begin{eqnarray*}
q^{n+1/2} & = & q^n + p^n \Delta t/2, \\
p^{n+1} & = & p^n - \Delta t \nabla V(q^{n+1/2}), \\
q^{n+1} & = & q^n + p^{n+1} \Delta t/2,
\end{eqnarray*}
 A practical steplength for simulations of our system with the St\"ormer-Verlet method is $\Delta t=0.01$.  This choice of $\Delta t$ is close to the largest for which the system can be integrated without an explosive instability in energy on the interval $[0,1000]$ for the initial conditions we consider. 

In the introduction we mentioned that researchers only consider low dimensional information from a molecular dynamics simulations. 
For our numerical experiments we will consider the configuration over time of the first particle: $q_1(t) \in \mathbb{T}^2, t \in [0,T]$. 
For the purposes of our experiments it helps to view the motion of the particle as occurring  in $\mathbb{R}^2$ and starting at the origin.
To this end, for the exact trajectory we define
\[
Q(t) = \int_0^t p_1(s) \, ds,
\]
and let $Q_x(t)$ and $Q_y(t)$ denote the respective $x$ and $y$ coordinates.  We have  that $Q_x(0)=Q_y(0)=0$ and 
if we let $t$ vary in $[0,T]$ then $Q \in (C[0,T])^2$. 
Similarly, for the numerical trajectory for each $\Delta t$ we define
\[
Q^n = \sum_{i=0}^{n-1} p^i.
\]
We then define $Q_{\Delta t}$ to be  the linear interpolation the the $Q^n$ at the times $n \Delta t$ so that $Q_{\Delta t} \in (C[0,T])^2$.

\begin{figure} 
\epsfig{file=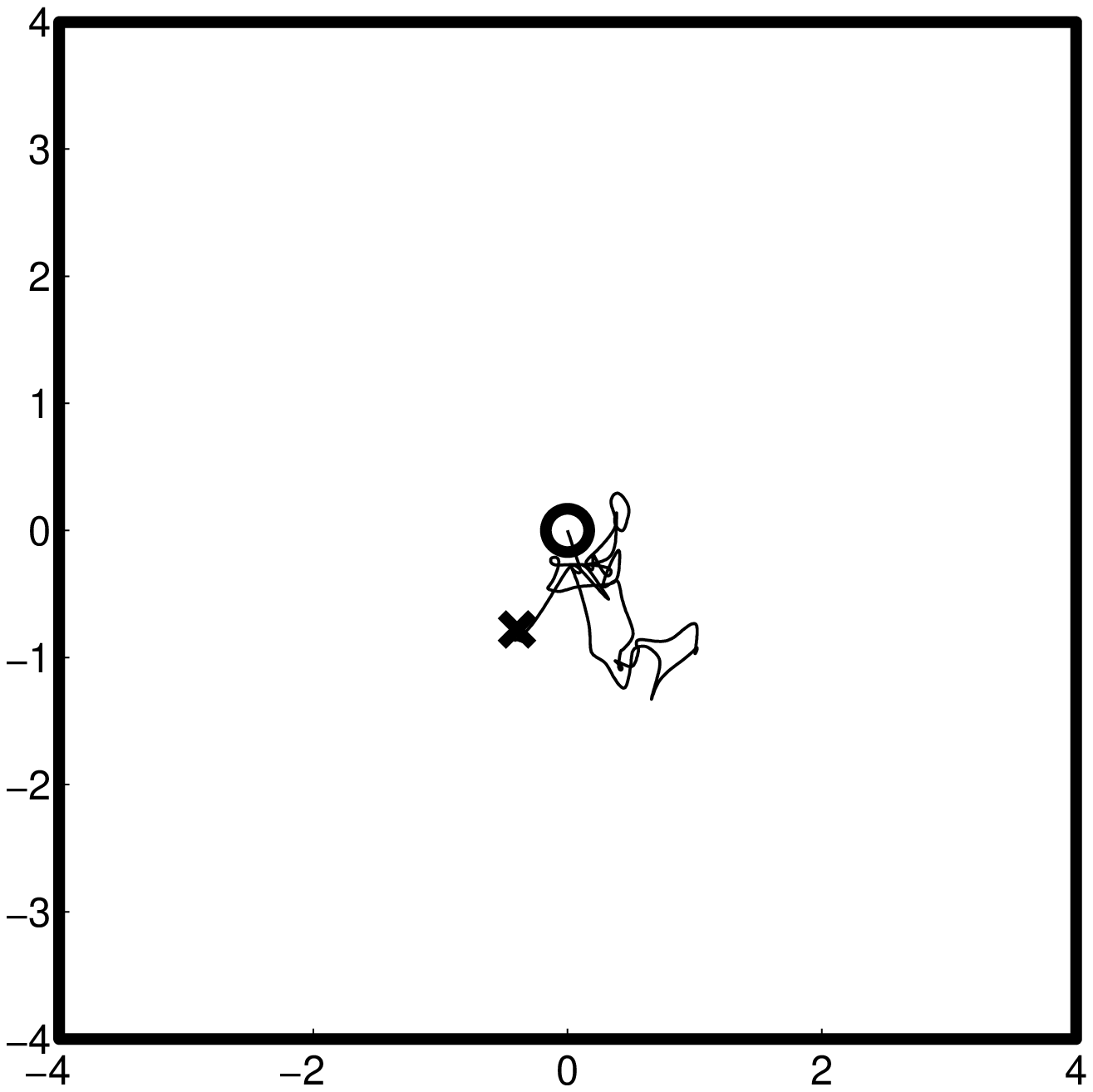,width=1.65in}
\epsfig{file=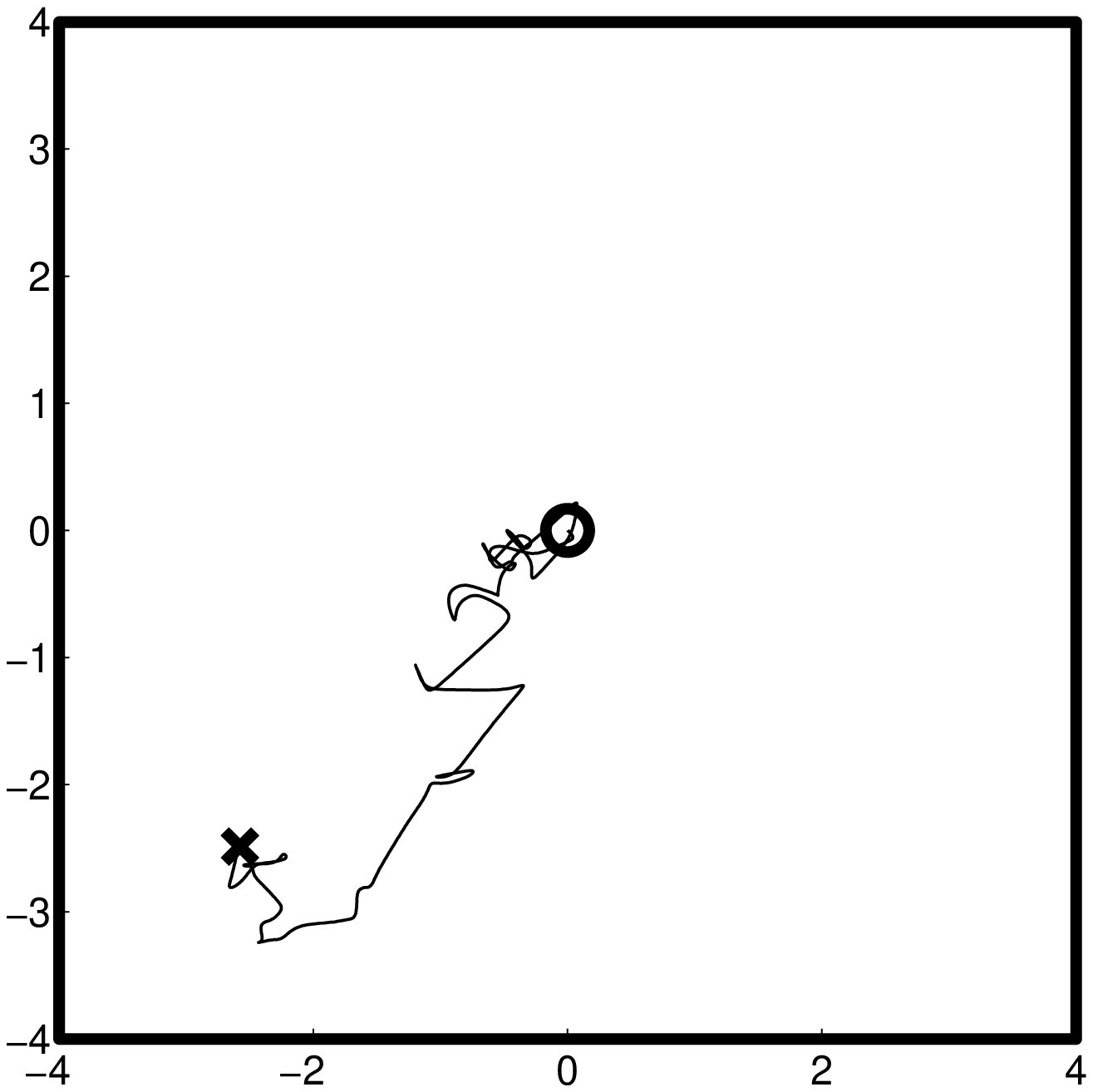,width=1.65in}
\epsfig{file=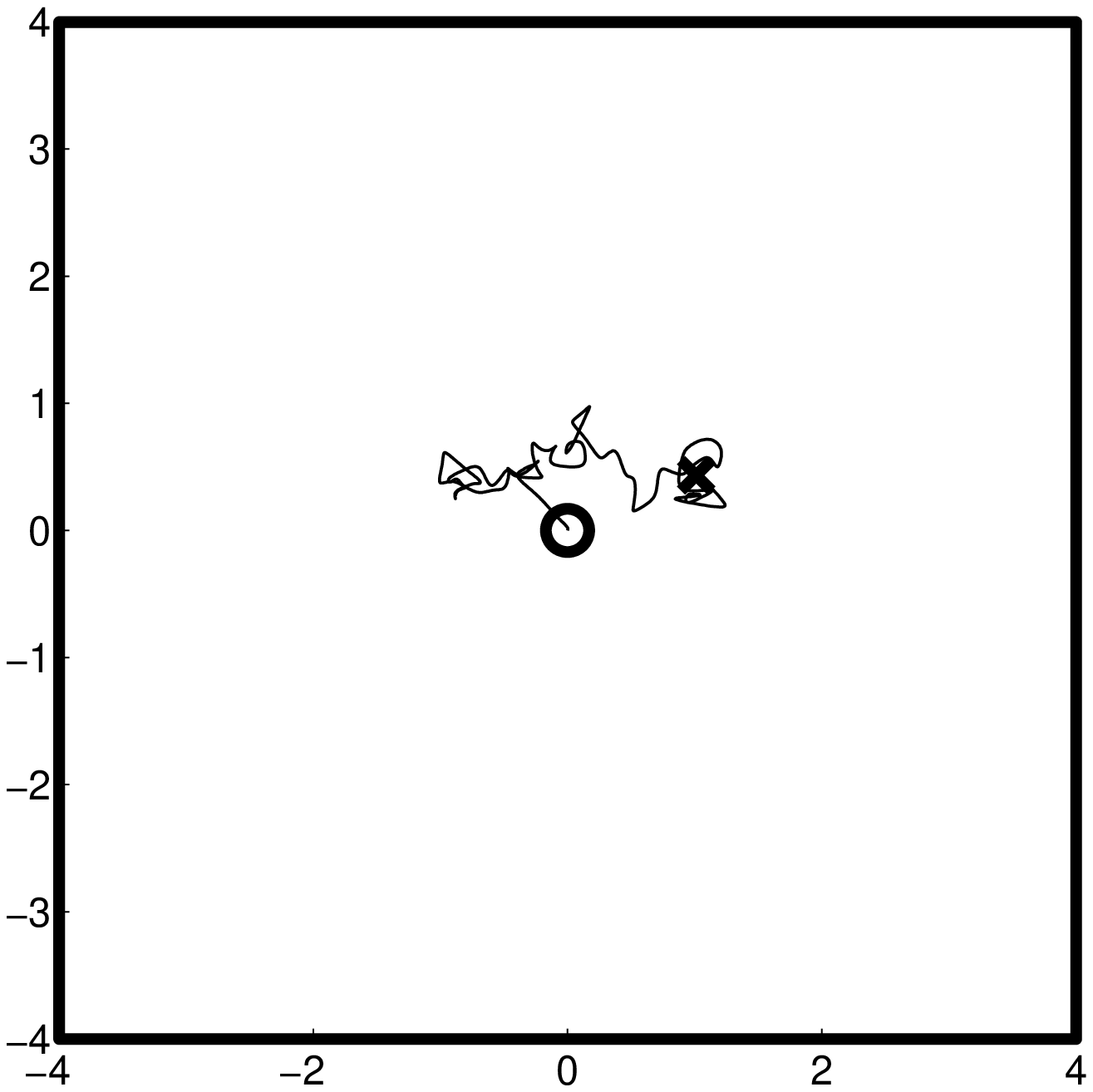,width=1.65in}
\caption{\label{fig:trajdemo}  First set of experiments.  Plots of realizations of  $Q_{\Delta t,x}$ vs.\ $Q_{\Delta t,y}$ for three different initial conditions over the time interval $[0,20]$ with $\Delta t= 0.01$.  The initial position is designated by a circle and the final position by an X.}
\end{figure}

 Our first set of numerical experiments demonstrates the qualitative features of the trajectories $Q$ over the time interval $[0,20]$.  We select three random initial conditions from the distribution given by (\ref{eq:distrib}) and plot in Figure~\ref{fig:trajdemo} the resulting numerical trajectories when $\Delta t=0.01$.   At the scale shown here there was not a noticeable qualitative difference between these plots and the similar plots generated with a smaller $\Delta t$.
 The motion is highly irregular, looking somewhat like Brownian motion in $\mathbb{R}^2$.   However, unlike Brownian motion, the exact trajectories $Q$ are infinitely smooth.  The interpolated numerical approximations $Q_{\Delta t}$ are piecewise linear.
  
\begin{figure} 
\epsfig{file=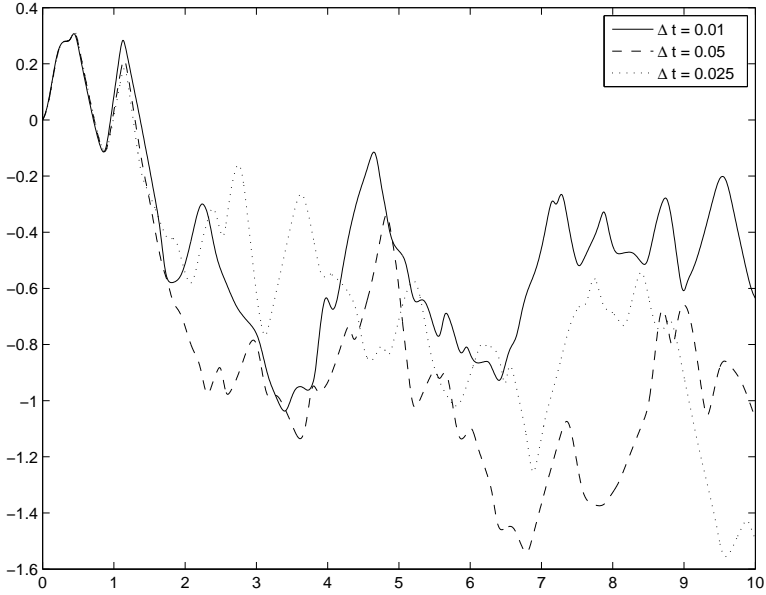,width=4in}
\caption{\label{fig:convergtraj} Second set of experiments.
Computed $Q_{\Delta t}(t)$
  versus $t$ for fixed initial conditions for a  range of $\Delta
  t$.}   
\end{figure}

Our second set of experiments demonstrates that individual trajectories computed using the timestep $\Delta t= 0.01$ are not accurate over time scales of interest.
 We randomly generate one initial condition according to the canonical distribution and then simulate over $[0,10]$ for $\Delta t = 0.01, 0.005, 0.0025$.   In Figure~\ref{fig:convergtraj} we plot $Q_{\Delta t,x}(t)$ versus $t$ for each of these steplengths.
 If the trajectory computed with steplength $\Delta t=0.01$ is
accurate over the time interval $[0,10]$, we expect that reducing the
time step by a factor of two would not yield a significantly
different curve.   However, we see that the two curves for $\Delta
t=0.01$ and $\Delta t=0.005$ very quickly diverge.   Moreover, we see that the trajectory with $\Delta t=0.005$ is not accurate either, since it diverges quickly from the trajectory with timestep $\Delta t=0.0025$.  Obtaining an accurate trajectory over the interval $[0,10]$ and certainly over $[0,100]$ would require $\Delta t$ to be considerably smaller than what is used in practice.  This same convergence behaviour is observed for all initial conditions selected according to the canonical distribution.

Our third set of experiments shows that despite the inaccuracy of individual simulations of the system, the statistical features of numerical trajectories appear to be highly accurate even for $\Delta t=0.01$.  Again we consider the trajectory of a single particle for  initial conditions drawn from the distribution defined by (\ref{eq:distrib}).   For each randomly generated initial condition we compute the value of a collection of functionals of the numerical trajectories over the time intervals $[0,100]$ and $[0,1000]$.  We plot these values in histograms in order to observe the distribution of each functional.  We do this for each of $\Delta t=0.01, 0.005, 0.0025$ and for five functionals of the trajectory.
The functionals $F_i \colon (\mathbb{R}[0,T])^2 \rightarrow \mathbb{R}$ we consider are 
\begin{eqnarray*}
F_1(Q) & = & Q_x(T), \\
F_2(Q) & = & \frac{1}{T} \int_0^T Q_x(t) \sin( 2 \pi t / T)\, dt, \\
F_3(Q) & = & \max_{t \in [0,T]} \| Q(t) \|, \\
F_4(Q) & = & \min \{ t  \in [0,T] \colon \| Q(t) \| \geq 1 \}, \\
F_5(Q) & = & \frac{(Q(T) - Q(T-\tau) )^{\mathrm{T}} ( Q(T-\tau) - Q(T-2\tau))}
{ \| Q(T) - Q(T-\tau) \| \, \| Q(T-\tau) - Q(T-2\tau) \| }, \ \ \ \tau = 0.1.
\end{eqnarray*}
$F_1$ is simply the $x$-position of the particle at time $T$.  $F_2$ is the average of  $Q_x(t) \sin(2\pi t/T)$ over $[0,T]$.  $F_3$ is the maximum distance from its initial condition that the  particle attains on $[0,T]$.  $F_4$ is the first time at which the particle leaves a ball of radius 1 centred at its initial condition.  Its value was set to $T$ if the particle did not leave within $[0,T]$.  $F_5$ is the cosine of the angle between two adjacent increments of $Q$ just before time $T$.

\begin{figure}
\begin{center}
\epsfig{file=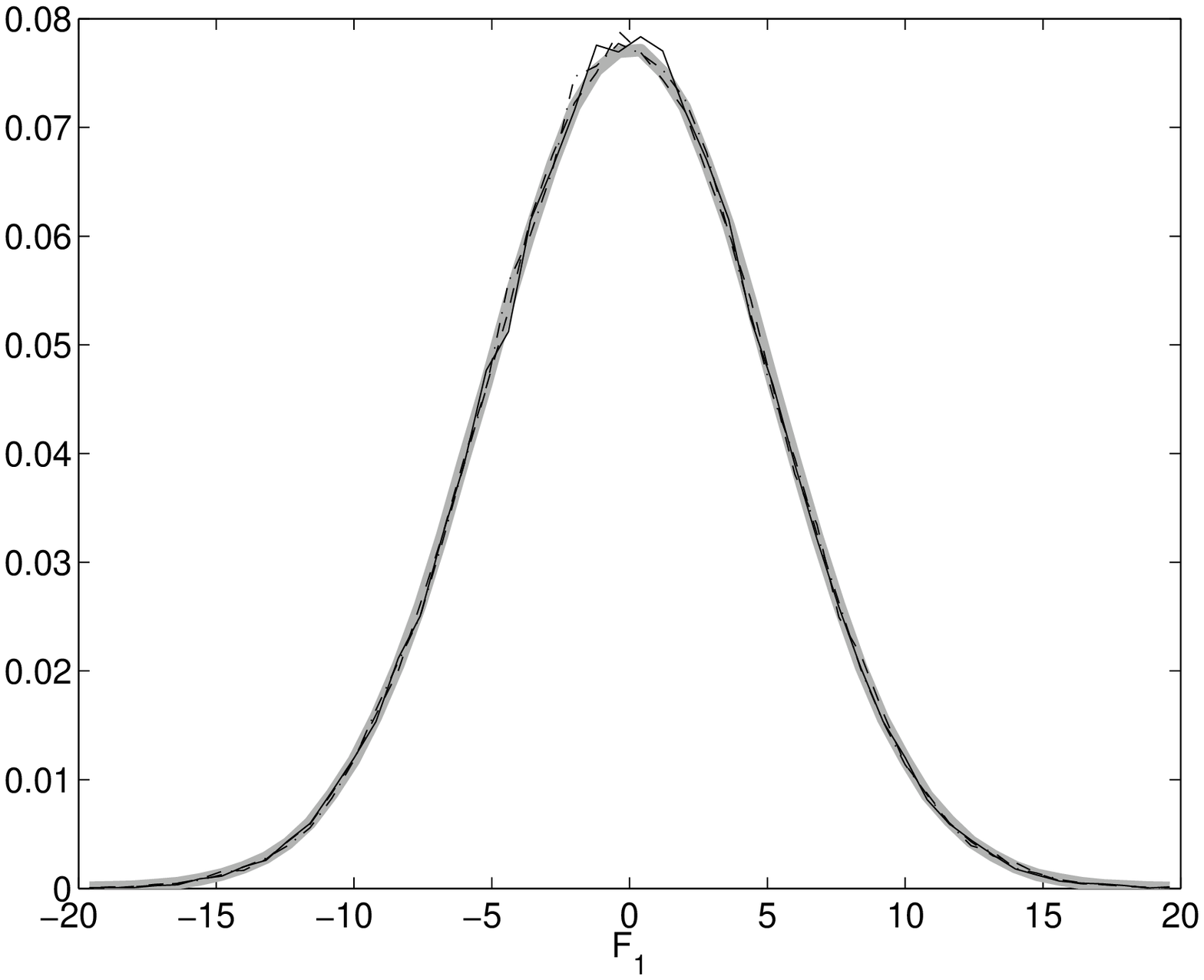,width=2.5in}
\epsfig{file=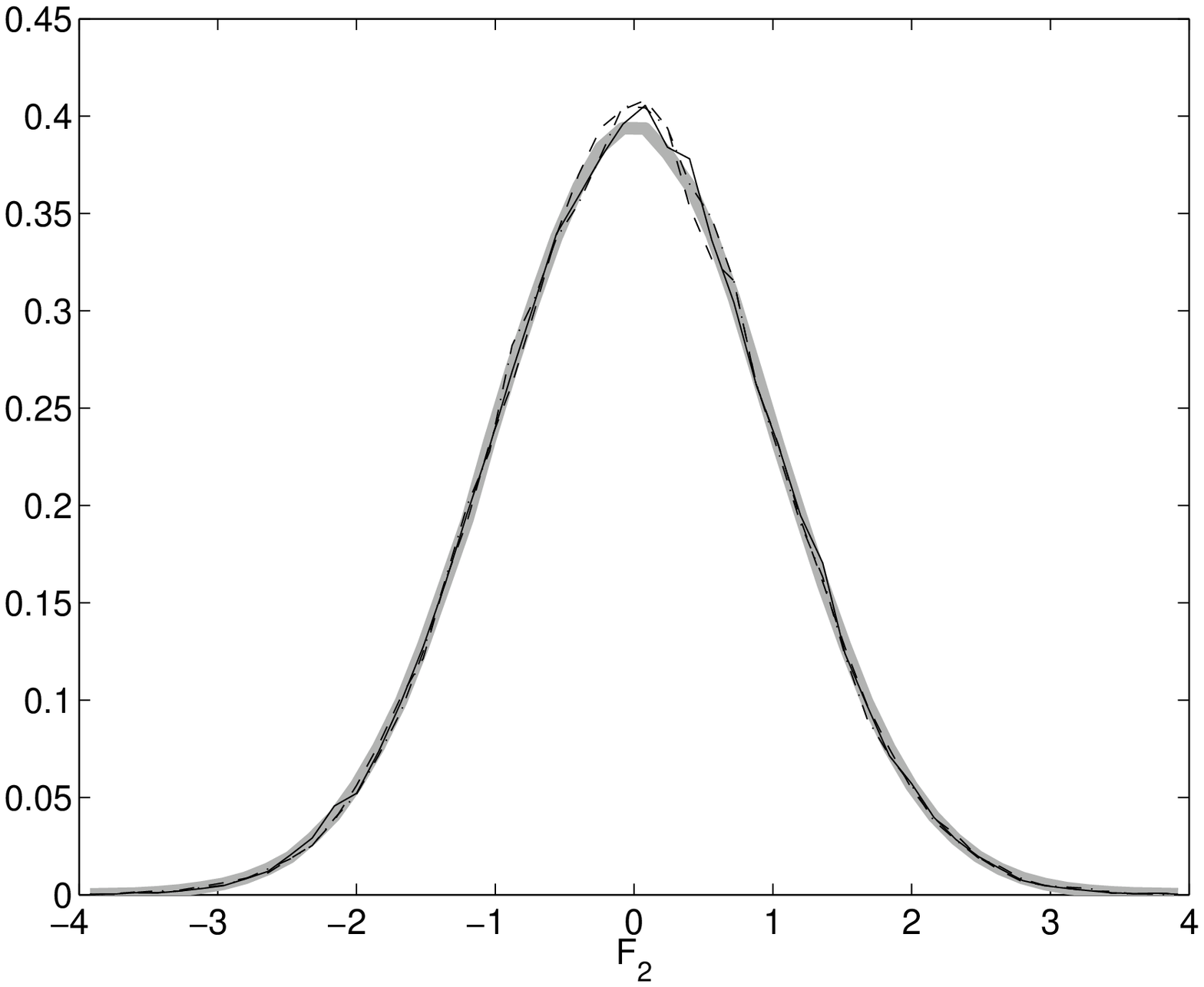,width=2.5in}\\
\epsfig{file=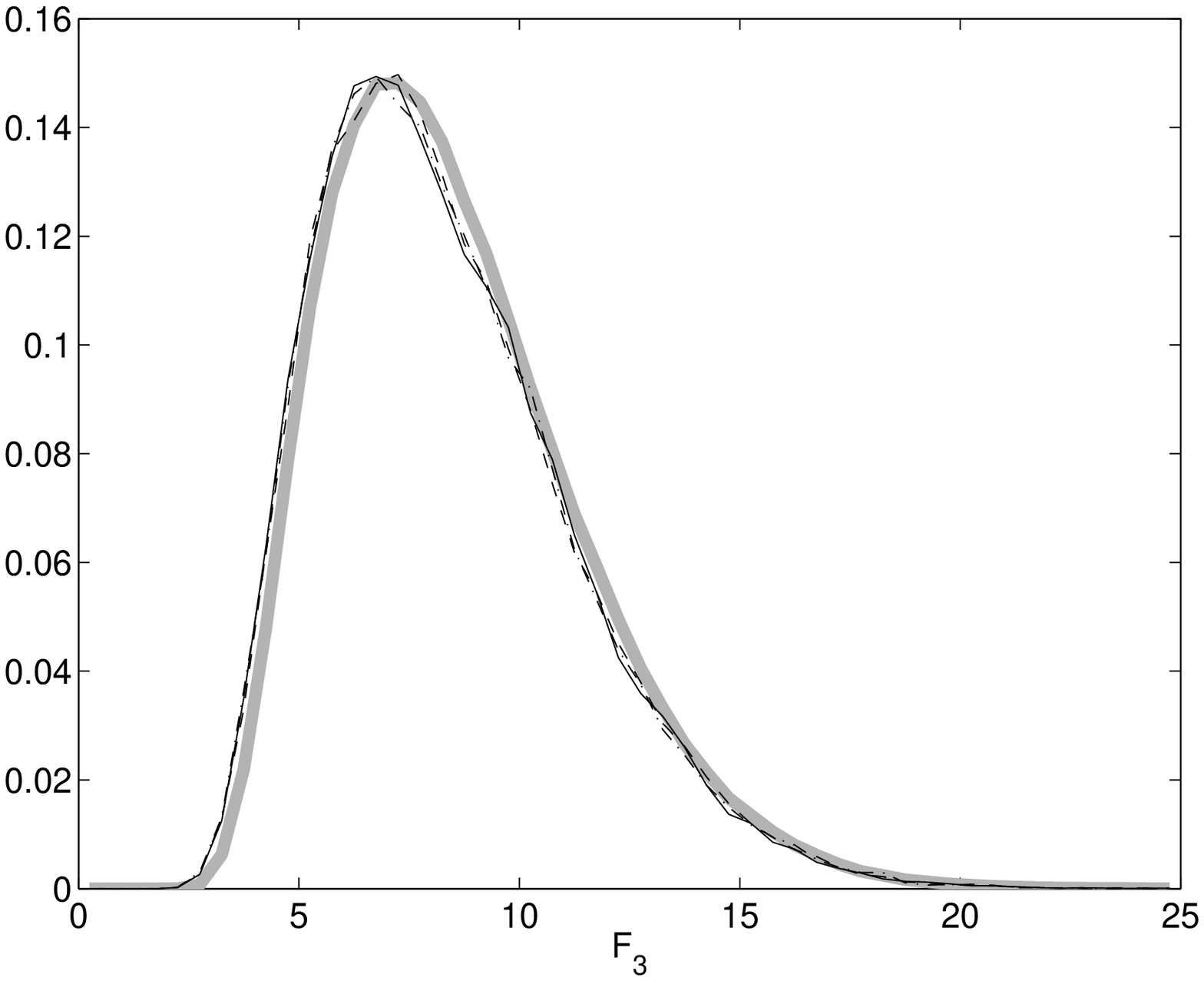,width=2.5in}
\epsfig{file=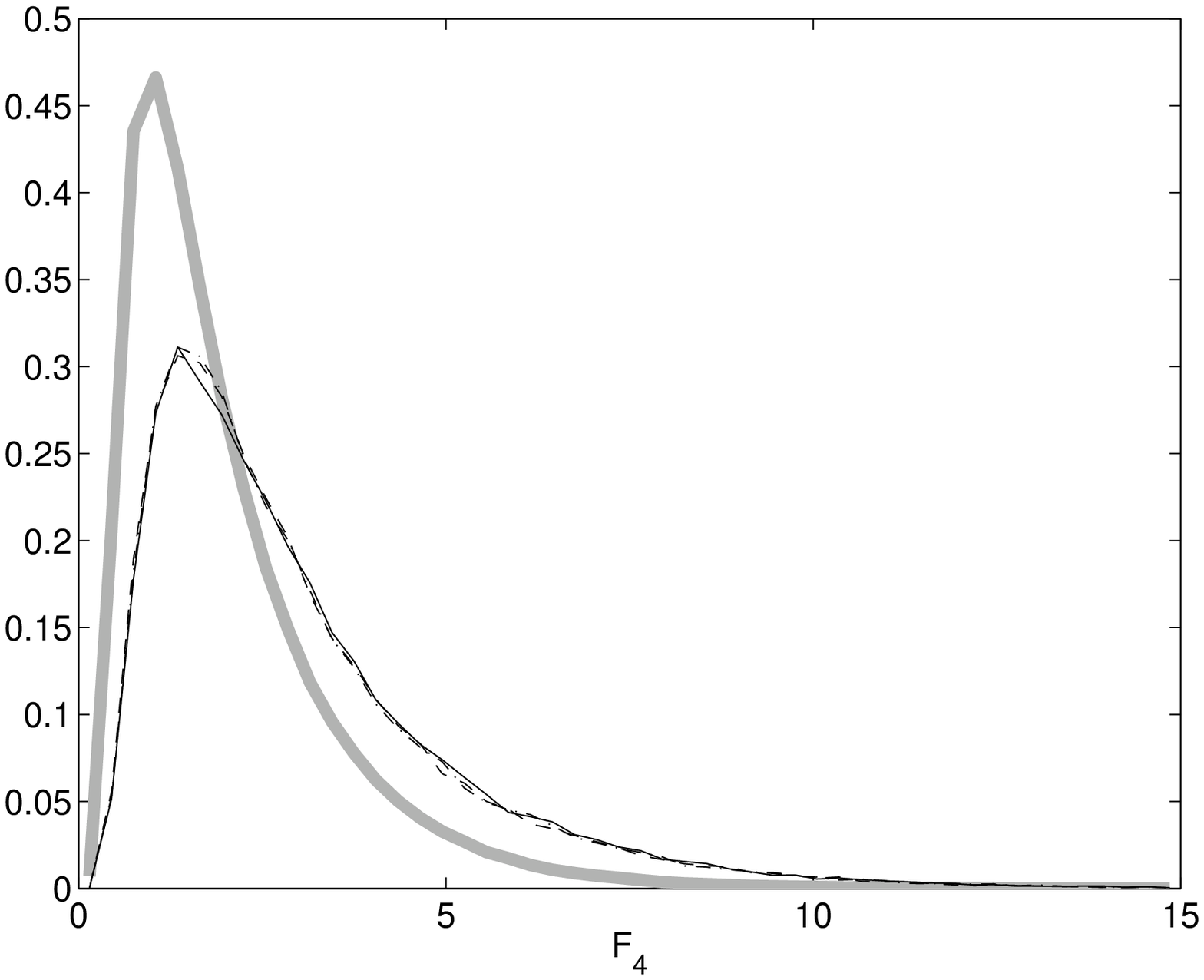,width=2.5in}\\
\epsfig{file=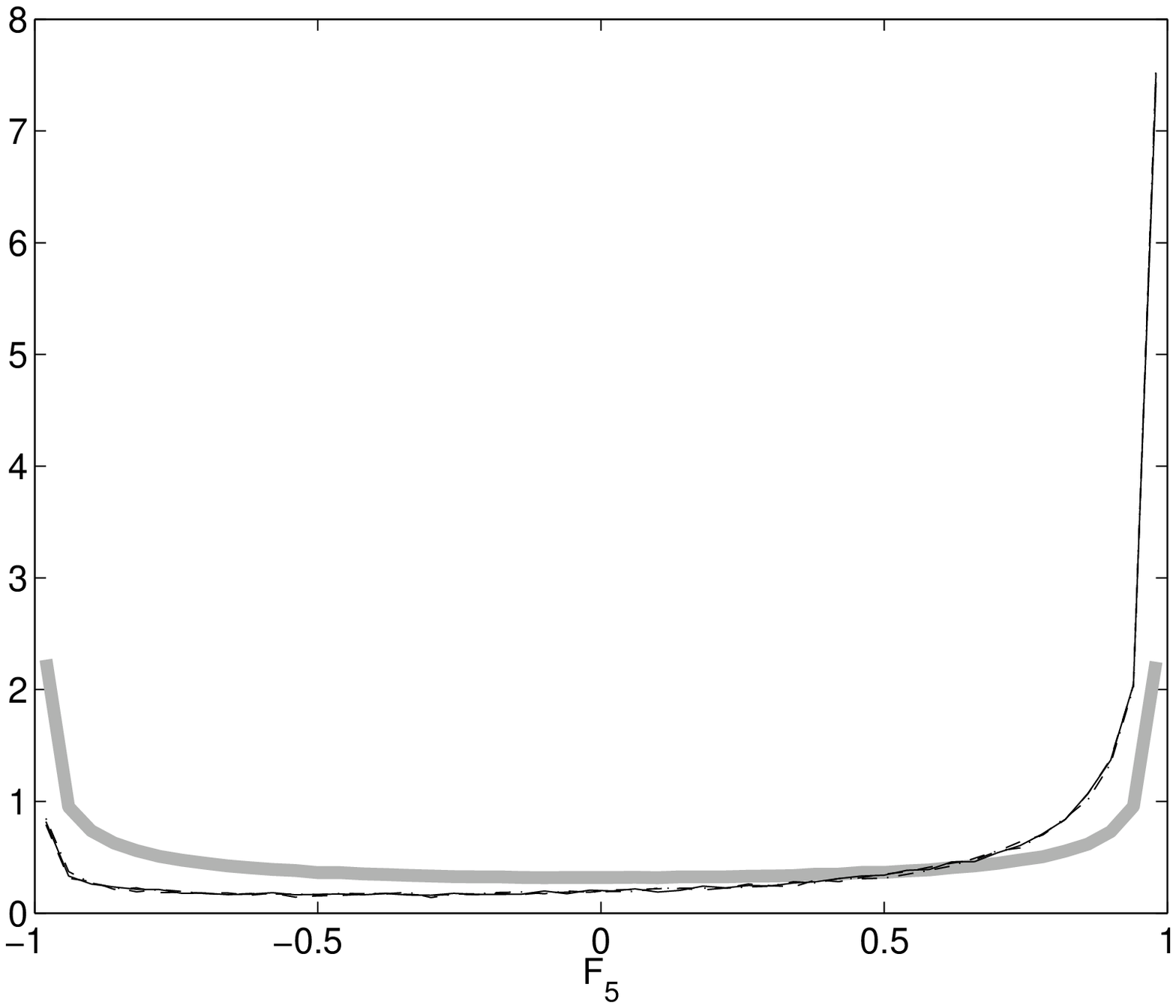,width=2.5in}
\end{center}
\caption{ \label{fig:hists} Third set of experiments for the time interval $[0,100]$.  Histograms of $F_i(Q_{\Delta t})$ for $i=1,\ldots,5$.  Each plot shows the result for the numerical trajectory with $\Delta t=0.01$ (solid), $\Delta t=0.005$ (dashed), $\Delta t=0.0025$ (dash-dot), as well as for Brownian motion (grey). } 
\end{figure}

\begin{figure}
\begin{center}
\epsfig{file=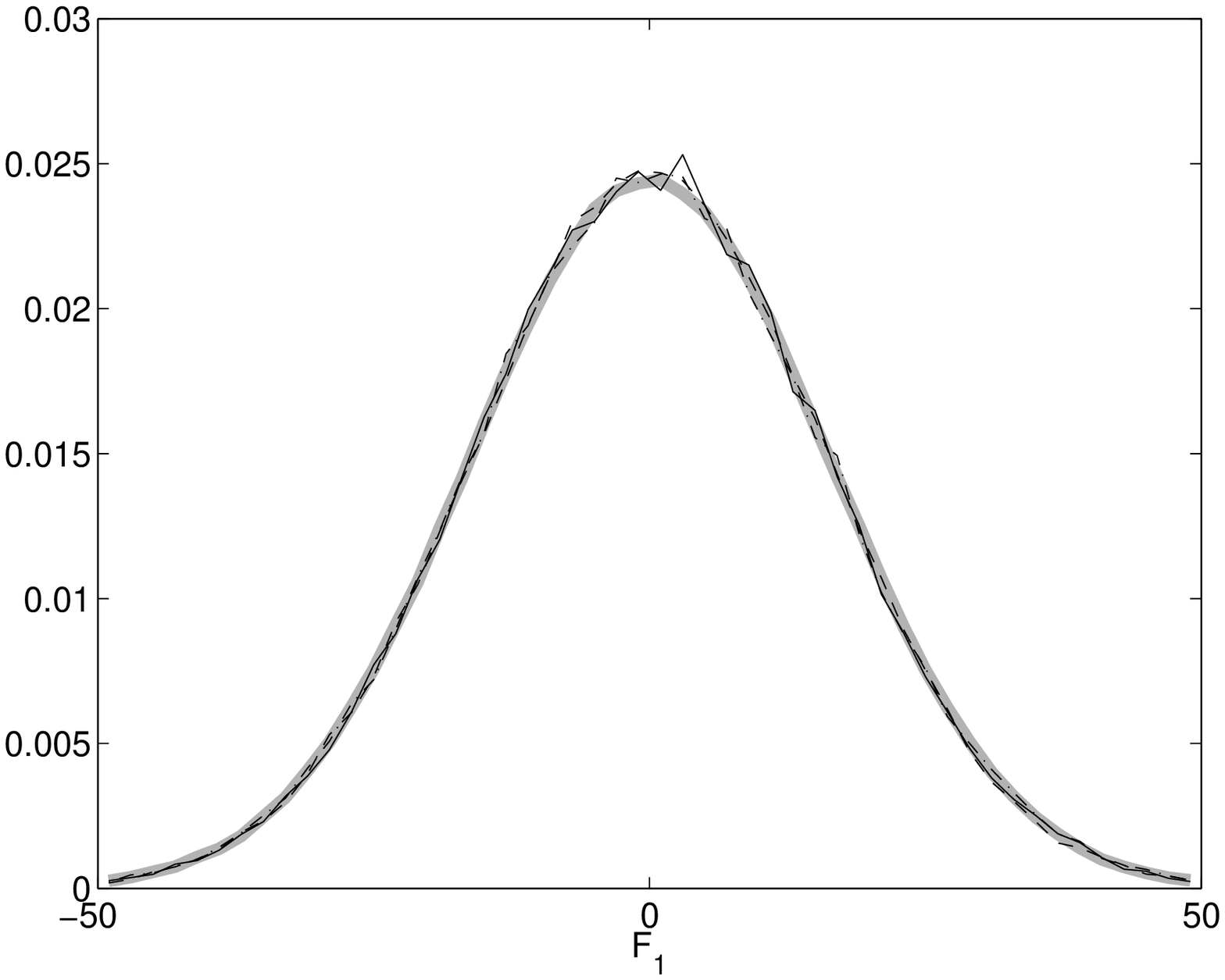,width=2.5in}
\epsfig{file=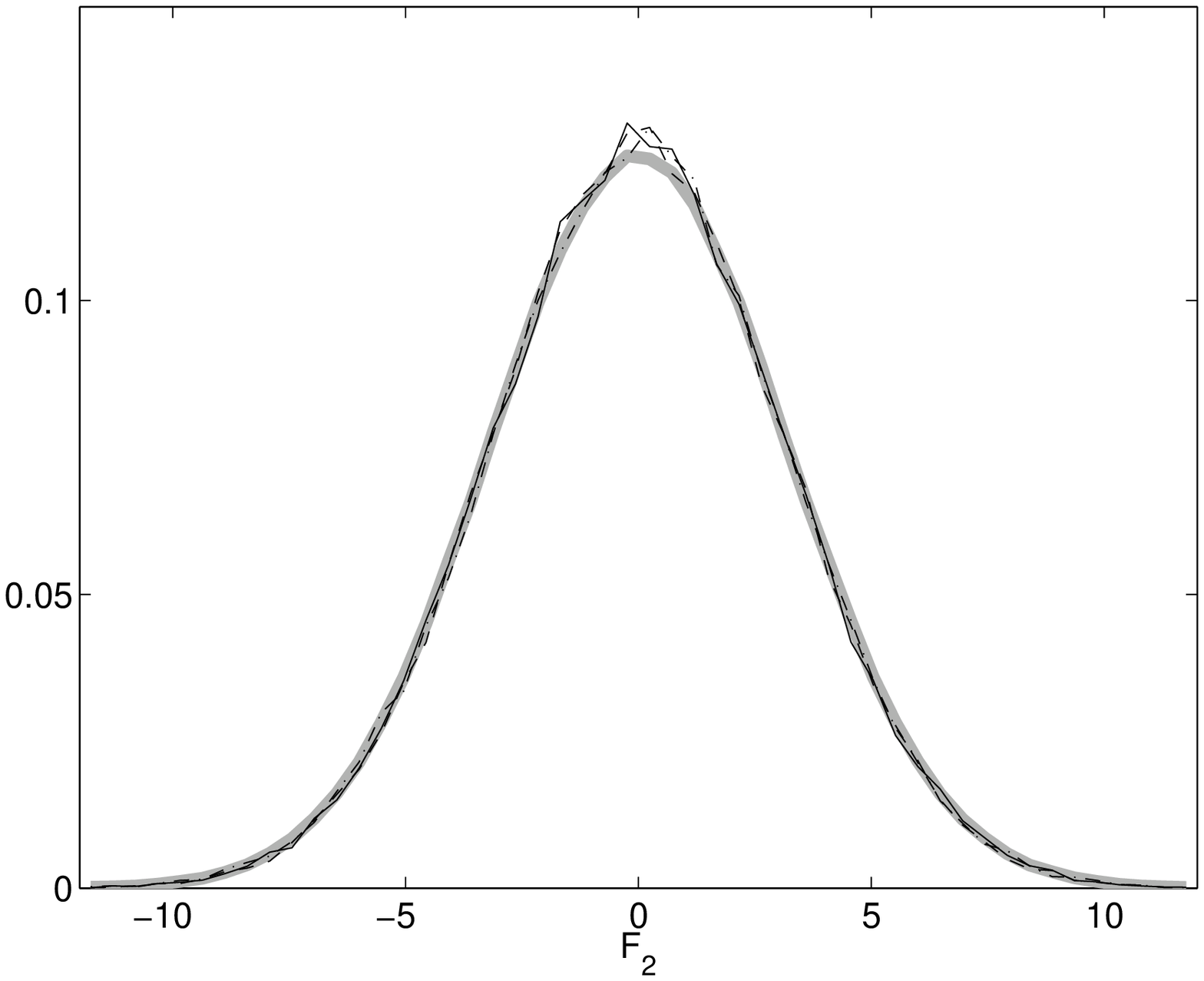,width=2.5in}\\
\epsfig{file=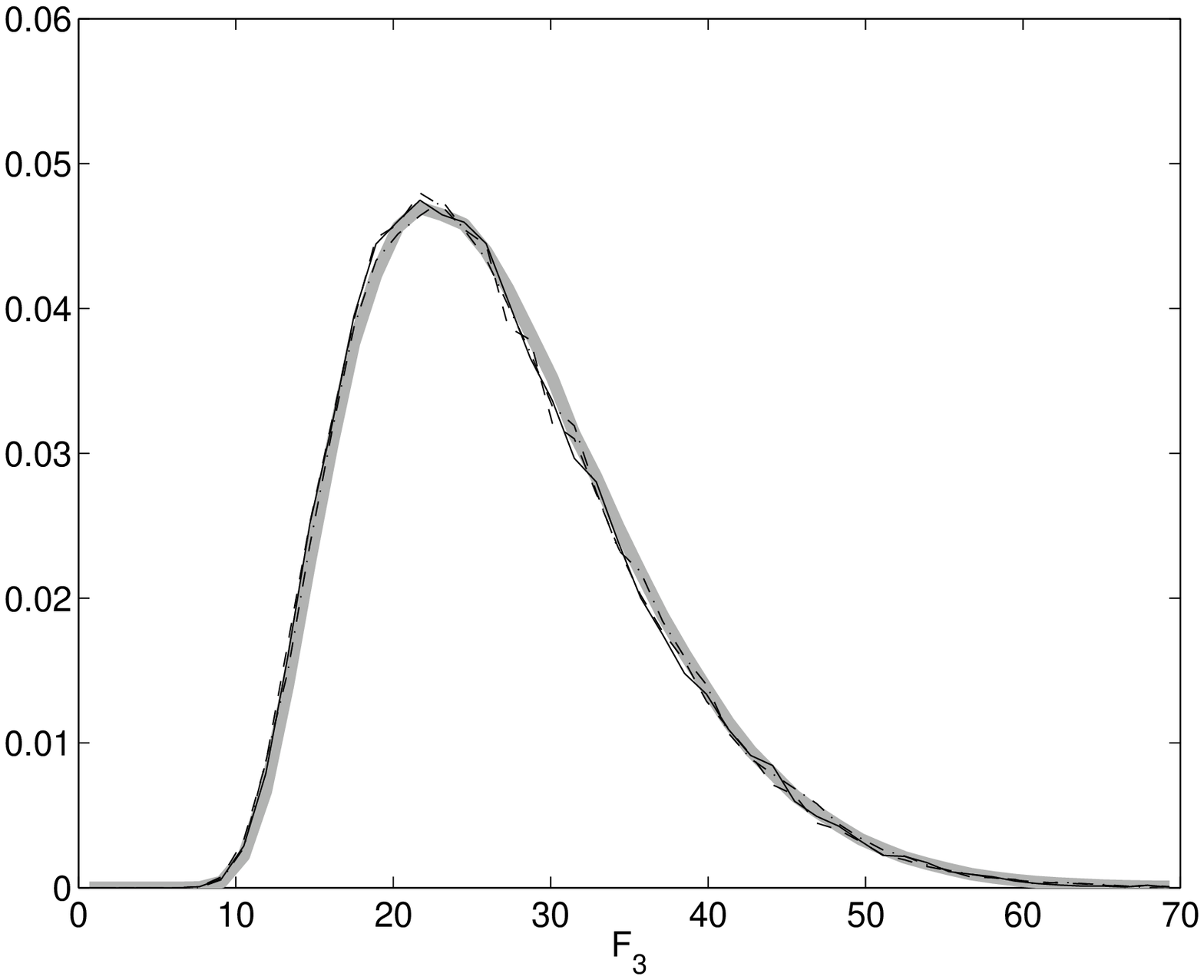,width=2.5in}
\epsfig{file=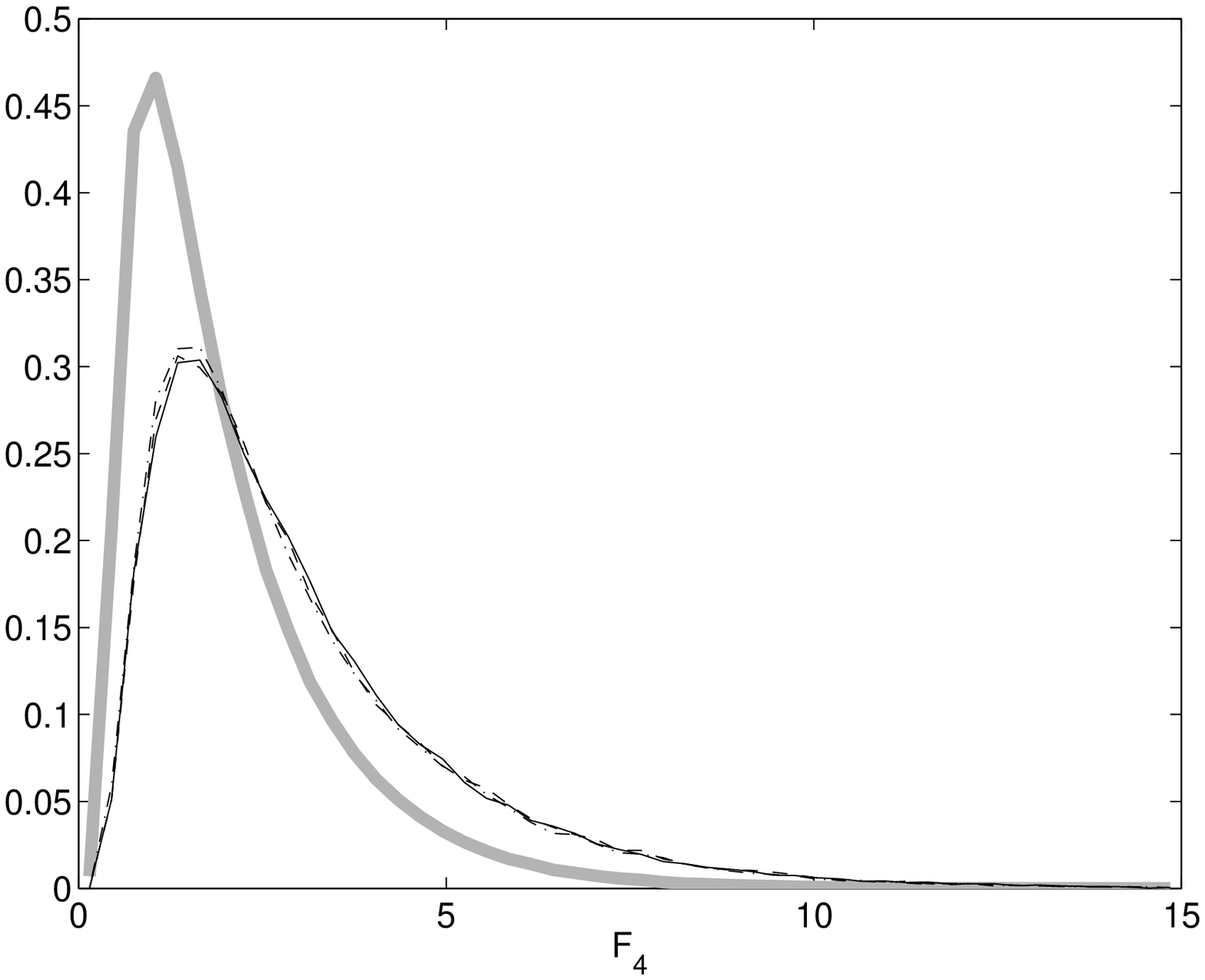,width=2.5in}\\
\epsfig{file=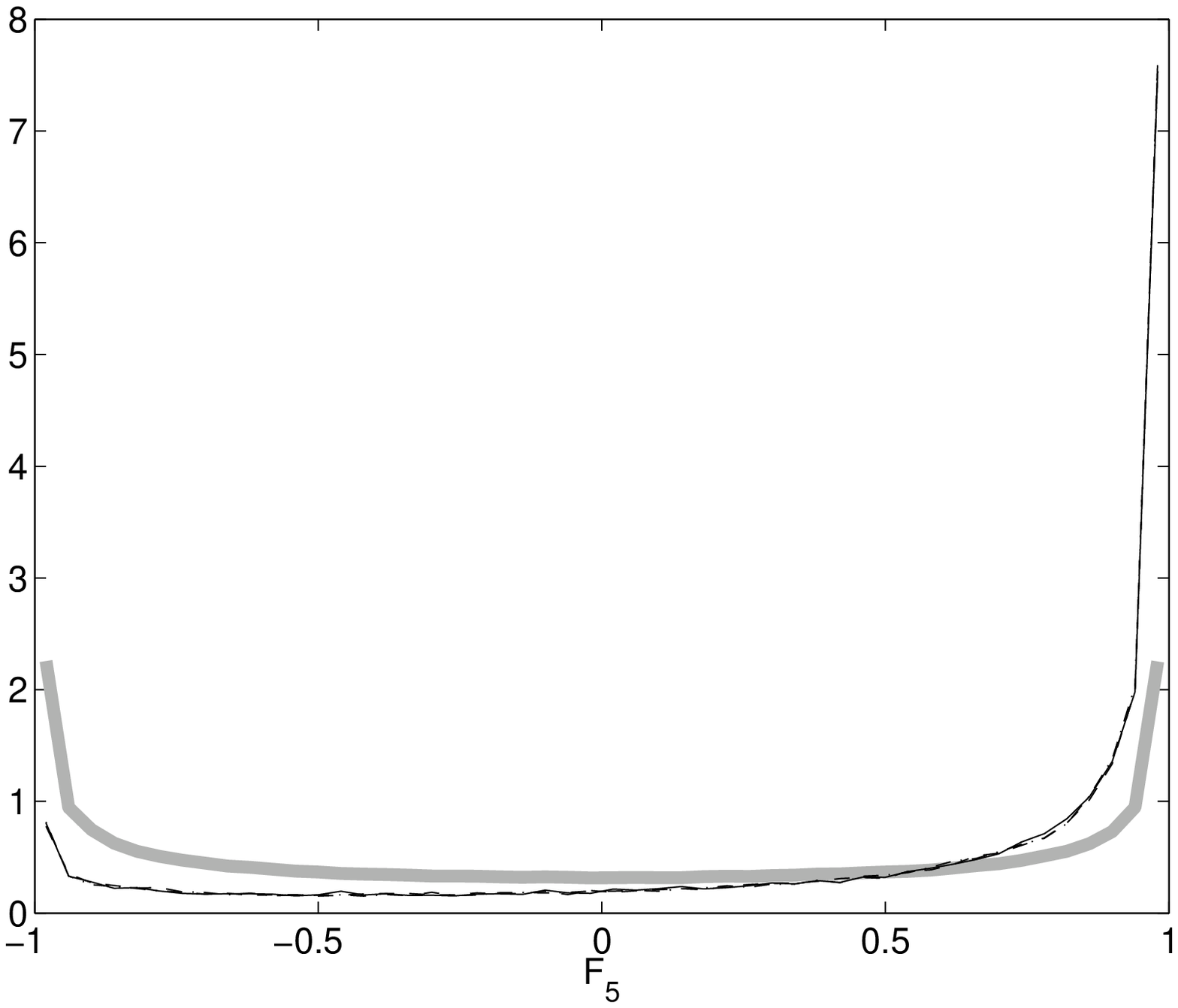,width=2.5in}
\end{center}
\caption{ \label{fig:hists2} Third set of experiments for the time interval $[0,1000]$.  Histograms of $F_i(Q_{\Delta t})$ for $i=1,\ldots,5$.  Each plot shows the result for the numerical trajectory with $\Delta t=0.01$ (solid), $\Delta t=0.005$ (dashed), $\Delta t=0.0025$ (dash-dot), as well as for Brownian motion (grey).} 
\end{figure}

In Figure~\ref{fig:hists} we show histograms of $F_i(Q_{\Delta t})$ for $i=1,\ldots,5$ with the different values of $\Delta t$ over the time interval $[0,100]$.  We also show the analogous histograms for two-dimensional Brownian motion $B(t)$ scaled so that $B_x(t)$ and $Q_{\Delta t,x}$ have the same variance.  We see that for all five functions the histogram generated does not appear to depend on the steplength used.  As well, we see that for some of the $F_i$ this matches closely the same histogram for Brownian motion, whereas for others it does not.  Note that we do not expect the histograms of $Q_{\Delta t}$ to converge to those for $B$: the exact trajectory $Q$ and Brownian motion $B$ do not have the same distribution.  These results are duplicated in Figure~\ref{fig:hists2} where we show the analogous plots for the time interval $[0,1000]$. 

\begin{figure}
\begin{center}
\epsfig{file=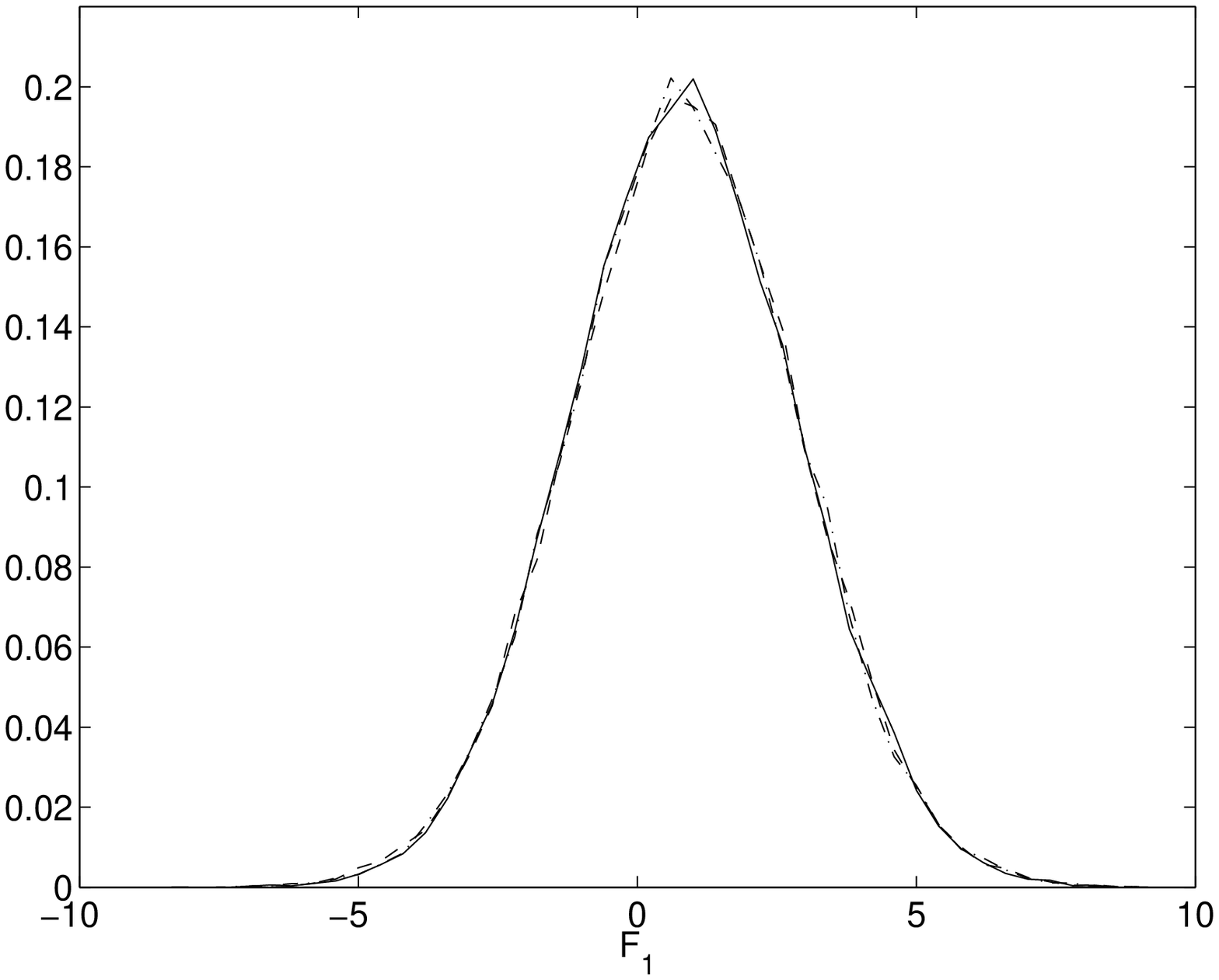,width=2.5in}
\epsfig{file=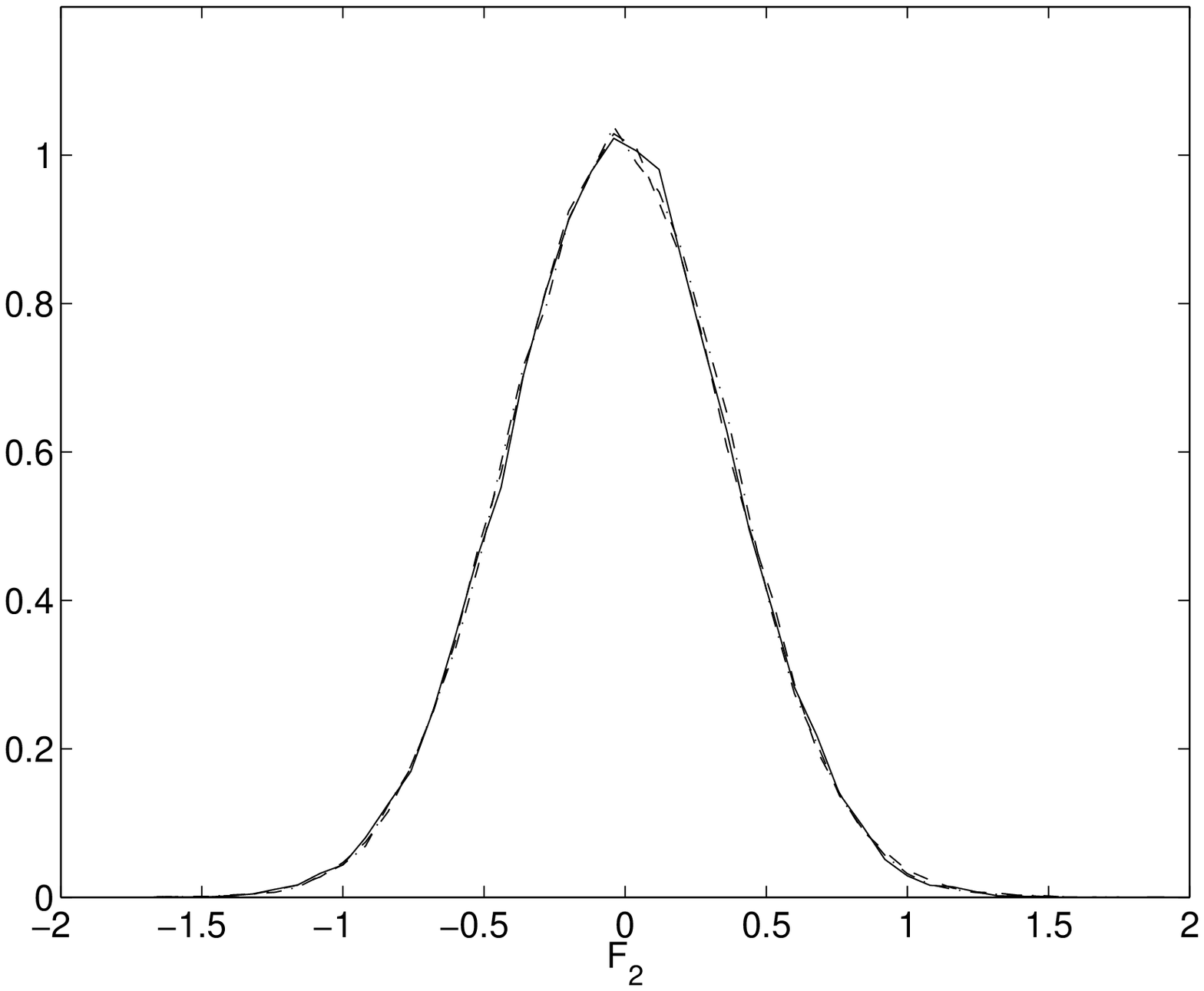,width=2.5in}\\
\epsfig{file=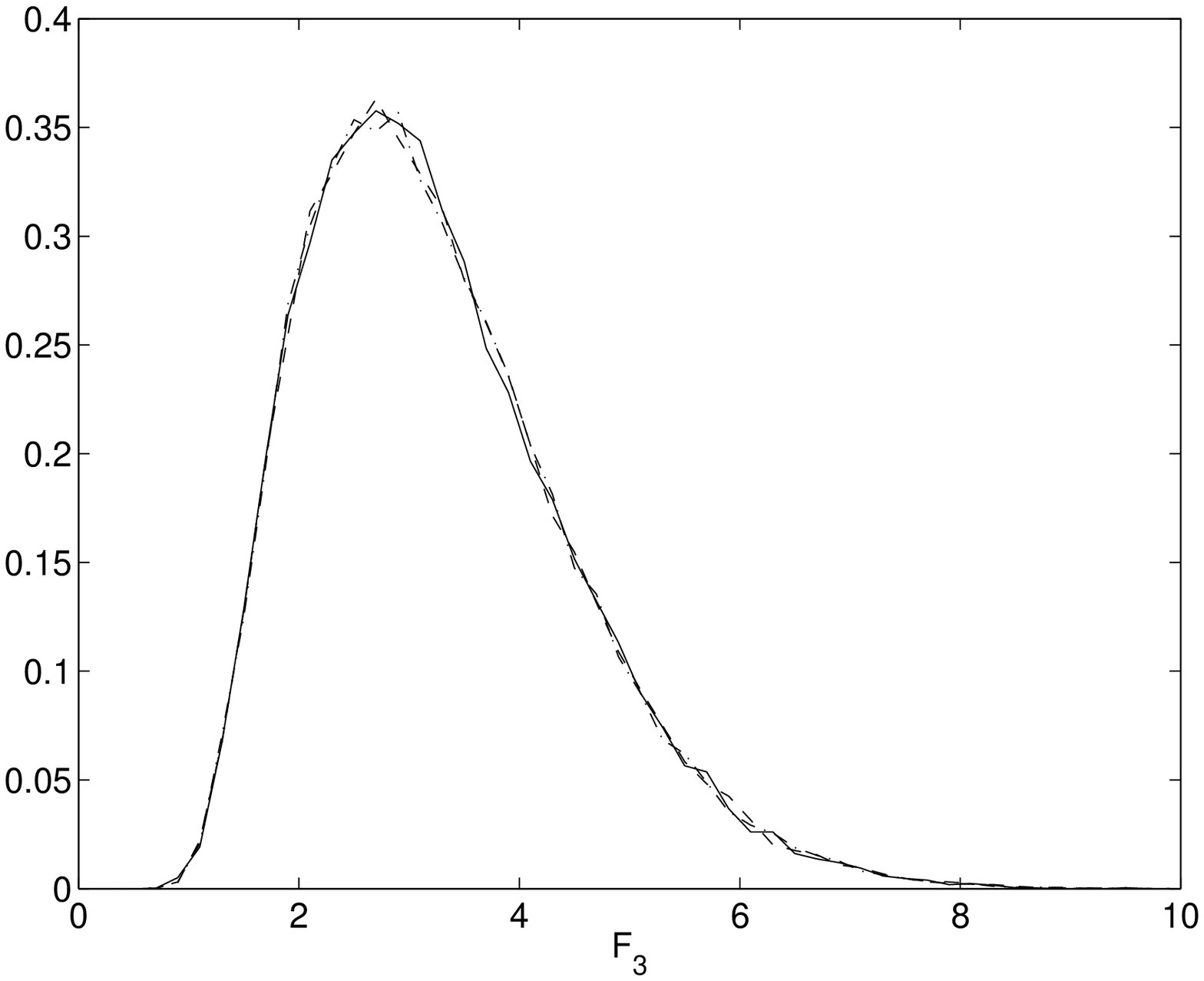,width=2.5in}
\epsfig{file=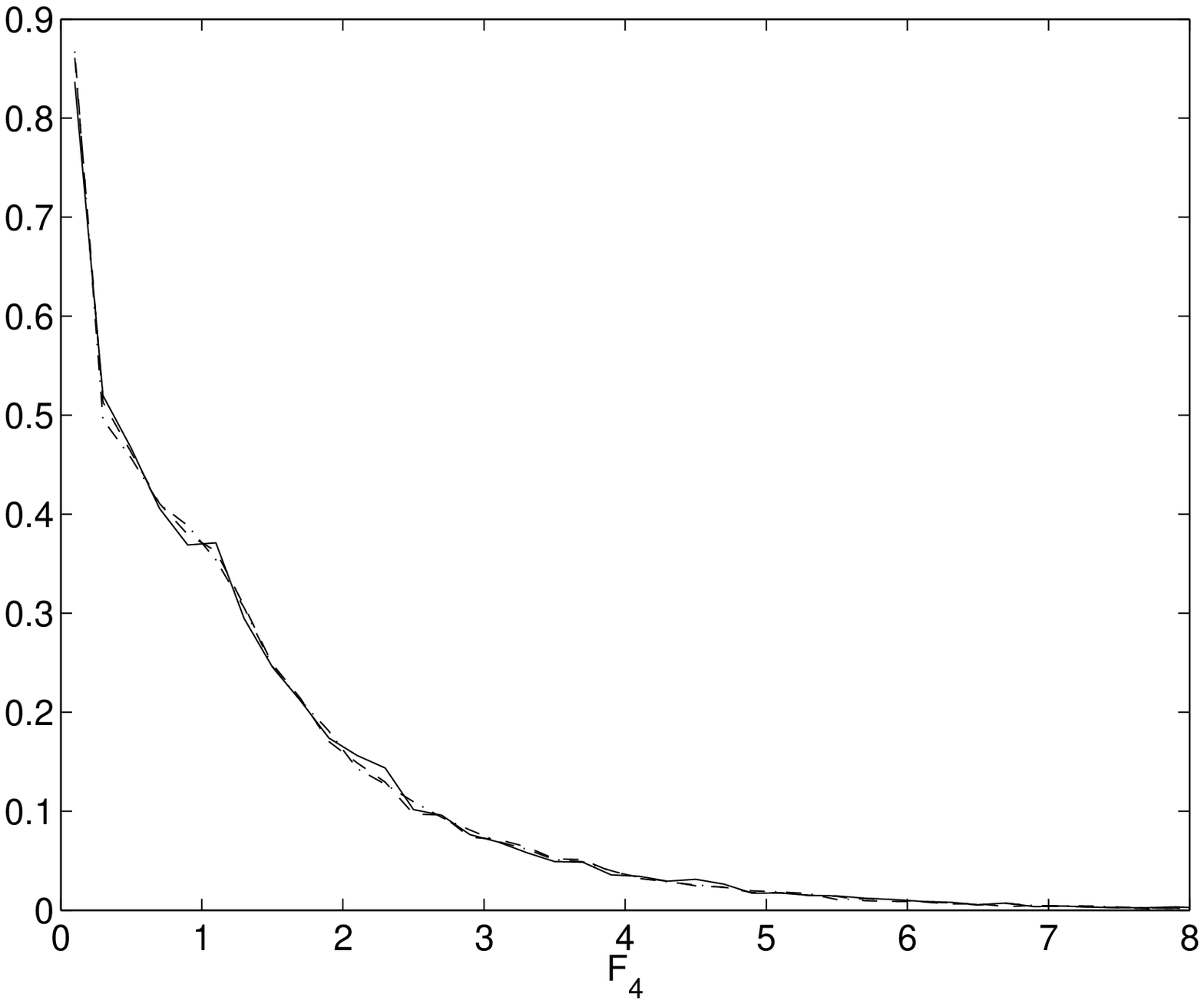,width=2.5in}\\
\epsfig{file=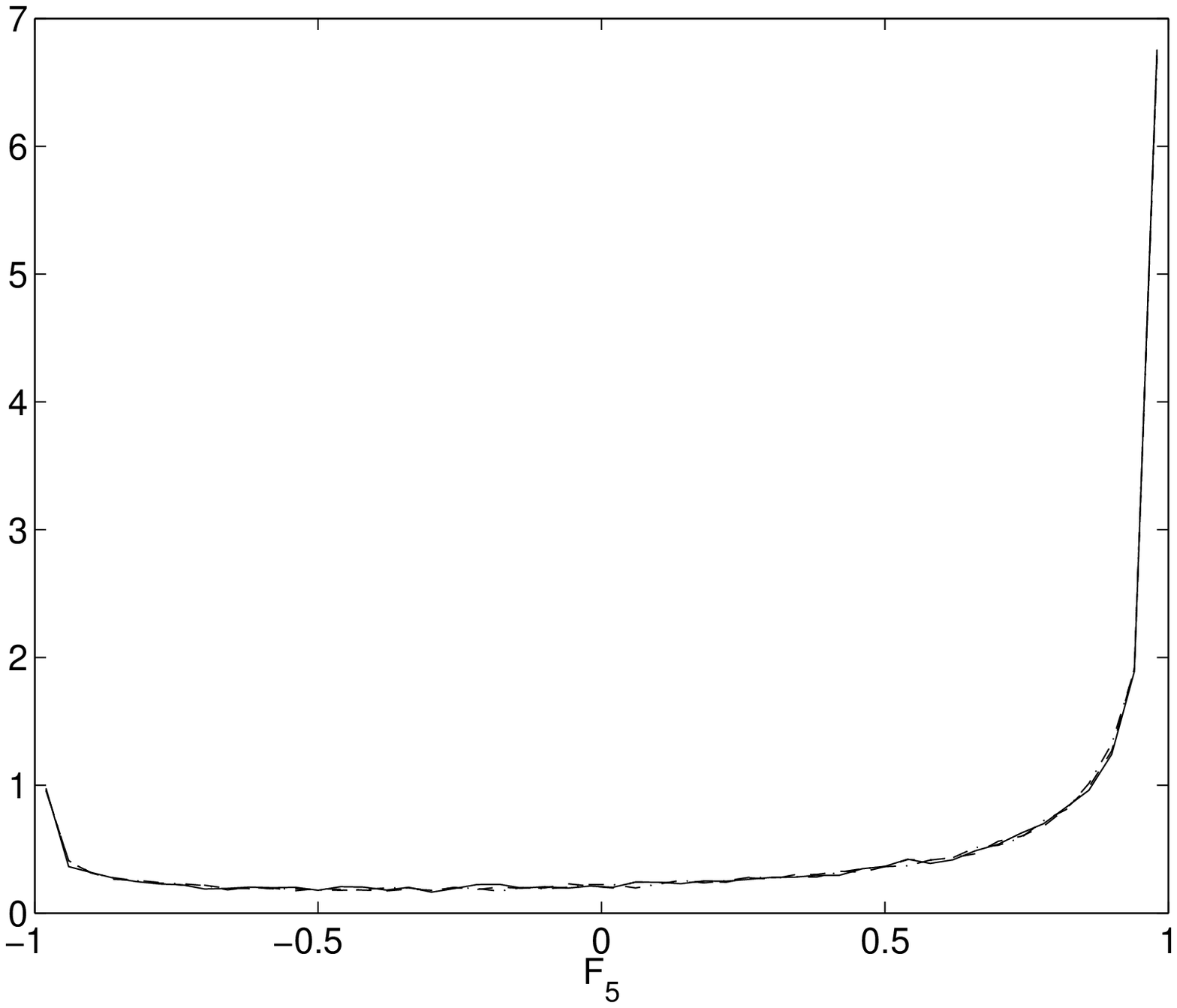,width=2.5in}
\end{center}
\caption{ \label{fig:noneq} Fourth set of experiments.   Using initial conditions drawn from a non-equilibrium distribution over  the time interval $[0,10]$.  Histograms of $F_i(Q_{\Delta t})$ for $i=1,\ldots,5$.  Each plot shows the result for the numerical trajectory with $\Delta t=0.01$ (solid), $\Delta t=0.005$ (dashed), $\Delta t=0.0025$ (dash-dot). } 
\end{figure}

Finally, in our fourth numerical experiment we repeat the previous experiment but we start with initial conditions that are drawn from a non-equilibrium distribution.  We randomly generate the initial conditions by first drawing from the equilibrium distribution (\ref{eq:distrib}) as before.  We then add 10 to the velocity in the $x$ direction of the first particle.
The typical trajectory arising from an initial condition selected in this way involves the first particle rapidly losing its excess energy through collisions with the other particles.
Within 10 time units the system is effectively indistinguishable from one started in the equilibrium state.
As before, we generate histograms of the functions $F_1, \ldots, F_5$, but now over the time interval  $[0,10]$ in order to highlight the effects of the nonequilibrium initial conditions.  Figure~\ref{fig:noneq} shows the results of these simulations.

For both the simulations from equilibrium and from non-equilibrium initial conditions
the histograms for all the functionals $F_1, \ldots, F_5$
are virtually identical for the different step-lengths, in contrast to the case where we examined single trajectories.  Any differences are well within the statistical error due to sampling only a finite number of trajectories.
This suggests that computed  distributions of the functionals with $\Delta t=0.01$ are fairly accurate for the distributions on the initial conditions we consider.  The non-rigorous argument for this belief is is as follows.  If the histograms were not accurate for the step length we were using, then reducing the steplength would cause the histogram to move significantly closer to the histogram for the exact solution.  Thus the histogram would not stay the same after halving the steplength.  The defect in this argument is that there could in theory be a broad range of values of $\Delta t$ for which the apparently same wrong histogram is computed.   An approximately correct histogram could be only observed for a $\Delta t$ much small than what we use.  
Despite this possibility, researchers generally trust histograms and other statistical information extracted from molecular dynamics trajectories \cite{frenkel}.

\section{Two Approaches} \label{sec:approaches}

Here we review two different proposals for the success of molecular dynamics: \emph{approximation in distribution} and \emph{weak shadowing}. 
For a distinct but related perspective, see \cite{skeelwhy}.

In the following $(C[0,T])^m$ denotes functions $x \colon [0,T] \rightarrow \mathbb{R}^m$.
 We use $X_0$ to denote a random initial condition in $\mathbb{R}^m$, $X \in (\mathbb{R}[0,T])^m$ to denote the random exact  trajectory of a dynamical system starting at $X_0$, and $X_{\Delta t} \in (C[0,T])^m$ the random approximate trajectory with the same initial condition.    We imagine the approximate trajectory to be generated by using a numerical method with steplength $\Delta t$ and then linearly interpolating the result.
 We measure the distance between two members $x,y$ of $(\mathbb{R}[0,T])^m$ by $d(x,y) = \| x-y\|_\infty = \sup_{t \in [0,T]} | x(t) - y(y) |$.
 Let $\Pi \colon \mathbb{R}^m \rightarrow \mathbb{R}^k$ be a map that extracts some low dimensional information from the system, so that the resulting low dimensional trajectories in $(C[0,T])^k$ are $\Pi(X)$ and $\Pi(X_{\Delta t})$.

\subsection{Approximation in Distribution} \label{subsec:distrib}

One explanation for the reliability of molecular dynamics is that if we let the initial condition of a simulation be random, then the distribution of the resulting numerical trajectory, seen as a random path in $(C[0,T])^m$, is close to the distribution of the actual trajectory.  We say that the trajectory is approximated in distribution.  This is also known as weak approximation.  Here we review some of the basic facts of approximation in distribution \cite{billingsley}.

Given a random variable $X$ taking values in $\mathbb{R}$, its distribution is the probability measure on $(\mathbb{R},\mathcal{B})$ defined by $\mu(A) = \mathbb{P}(X \in A)$ for $A \in \mathcal{B}$.  Here  $\mathcal{B}$ is the Borel subsets of $\mathbb{R}$.
  We say that two random variables $X$ and $Y$ have the same distribution if $\mathbb{P}(X \in A) = \mathbb{P}(Y \in A)$ for all $A \in \mathcal{B}$.  This is equivalent to $\mathbb{E} f(X) = \mathbb{E} f(Y)$ for all measurable $f$.
Note that two random variables need not be close on a realization-by-realization basis in order for their distributions to be identical.

The concept of distribution extends naturally to random vectors taking values in $\mathbb{R}^m$ and indeed to random elements of  any metric space as follows.
Consider a metric space $(\mathcal{S},d)$ with metric $d$ and let $\mathcal{B}$ be the  Borel subsets of $\mathcal{S}$ induced by $d$.  A random element of $\mathcal{S}$ is a measurable function $X: \Omega \rightarrow \mathcal{S}$ where $(\Omega,\mathcal{F},\mathbb{P})$ is some probability space. 
 The distribution of $X$ on $\mathcal{S}$ is the probability measure given by $\mu(A) = \mathbb{P}(X\in A)$.
As in the case of random variables, two random elements $X$ and $Y$ of $\mathcal{S}$ can have the same distribution without  $X(\omega)$ and $Y(\omega)$ being close for any fixed $\omega \in \Omega$.

Suppose we want to quantify how close the distributions of two random elements of a metric space are to each other.  A natural way to do this is to define a metric on the space of random elements of a metric space.    One popular choice is the Prokhorov metric, $\rho$, which we define here.  It has the property that if $\rho(X,Y)=0$ for random elements  $X$ and $Y$  if and only if $X$ and $Y$ have identical distributions.  

For a set $A \subset S$ and $\epsilon \geq 0$ we define $A^{\epsilon}$, the set of all points within distance $\epsilon$ of $A$ by 
\begin{equation}\label{eqn:epdef}
A^{\epsilon} = \{ x \in S \mid \inf_{y \in A}d(x,y) \leq \epsilon \} = \{ x \in S \mid d(x,A) \leq \epsilon\}.
\end{equation}
The Prokhorov metric is defined as follows.

\begin{definition} \label{def:prokhorov} \cite[p.\ 72]{billingsley} For random variables $X$ and $Y$ in $S$
\[
\rho(X,Y) : = \inf \left\{ \epsilon > 0 \mid \mathbb{P} (X \in A) \leq \mathbb{P}(Y \in A^{\epsilon}) + \epsilon \right\}.
\]
\end{definition}

If we identify random elements of $S$ that have the same distribution, then $\rho$ is a metric on the set of random elements \cite[p.\ 394]{dudley}.    If $(S,d)$ is separable (as are all examples in this paper) random elements $X_n$ converge in distribution to $X$ if and only if $\rho(X_n,X) \rightarrow 0$ \cite[p.\ 395]{dudley}.  Note that $\rho(X,Y) \leq 1$ always.

A straightforward way to measure how close $X$ is to $X_{\Delta t}$ in distribution is to consider $\rho(X,X_{\Delta t})$, where we view $X$ and $X_{\Delta t}$ as random elements of $(C[0,T])^m$.  
However, we generalize this idea by measuring how close the distributions  of low dimensional functions of the full trajectories are.  
We consider  the map $\Pi \colon \mathbb{R}^m \rightarrow \mathbb{R}^k$,
 and we measure $\rho( \Pi (X), \Pi (X_{\Delta t}))$, where
$\Pi(X)$ and $\Pi(X_{\Delta t})$ are random elements of $(C[0,T])^k$.  Choosing $k=m$ and $\Pi$ to be the identity gives the approximation in distribution of $X_{\Delta t}$ to $X$ so this is a generalization of the original idea.  
As an example,  suppose the differential equations defining $X$ describe the motion of a system of particles in $\mathbb{R}^2$, so that the dimension of the system is $4n$.  To study the position of one particle as a function of time, one could let $\Pi \colon \mathbb{R}^{4n} \rightarrow \mathbb{R}^2$ be the function that simply returns the position of the first particle in the system.  This choice in analogous to what we did for our model system in Section~\ref{sec:problem}.

With these definitions in mind, we are led to a 
quantification of our belief that $\Pi(X)$ and $\Pi(X_{\Delta t})$ are close in distribution: we conjecture that
\begin{equation}  \label{eqn:approxdist}
\rho( \Pi(X), \Pi(X_{\Delta t}) ) < \epsilon,
\end{equation}
for some small $\epsilon$.
For the conjecture to be applicable to molecular dynamics, we must be able to control $\epsilon$ and the length of the time interval $T$ in terms of $\Delta t$.
We conjecture that for all sufficiently small $\Delta t$ there are some $C,D,E,p>0$
\[
\rho( \Pi(X), \Pi(X_{\Delta t}) ) < C \Delta t^p,
\]
for $T <  D \exp( - E / \Delta t)$.

One of the consequences of approximation in distribution with respect to the Prokhorov metric is that it allows us to bound the error we make in computing the expectation of functionals of the paths.  Suppose that $G : (C[0,T])^k \rightarrow \mathbb{R}$ is a bounded Lipschitz continuous function of the paths.  Let the norm $\| \cdot \|_{BL}$ be defined on the set of all such $G$ by 
\[
\| G \|_{BL} := \sup_{x} |G(x)| + \sup_{x \neq y} \frac{| G(x)-G(y) |}{\|x-y\|_{\infty}},
\]
where the suprema are taken over all $x,y \in  (C[0,T])^k$ \cite[p.\ 390]{dudley}.
We can define a another metric on the space of random elements of $(C[0,T])^k$ by
\[
\beta(X,Y) := \sup_{\| G\|_{BL} \leq 1} | \mathbb{E} G(X) - \mathbb{E} G(Y) |,
\]
originally defined in \cite{fortet}.  For any two random elements $X$ and $Y$ of a metric space, we have \cite[p.\ 398]{dudley}
\begin{equation} \label{eqn:rhobeta}
\rho(X,Y) \leq \left( \frac{3}{2} \beta(X,Y) \right)^{1/2},
\end{equation}
 and \cite[p.\ 411]{dudley} 
\[
\beta(X,Y) \leq 2 \rho(X,Y).
\]
So (\ref{eqn:approxdist}) would  imply $\beta(\Pi(X),\Pi(X_{\Delta t}) < 2 \epsilon$, and so
\begin{equation} 
\label{eqn:goal}
| \mathbb{E} G(\Pi(X)) - \mathbb{E} G(\Pi(X_{\Delta t})) | < 2 \| G \|_{BL}\epsilon.
\end{equation}
Of course, most functions $G$ of interest  are not bounded, but similar results can be obtained for unbounded, locally Lipschitz $G$ in the case that $G(\Pi(X))$ and $G(\Pi(X_{\Delta t}))$ have finite moments; see  \cite{tuppersde} for an example.

\subsection{Weak Shadowing} \label{subsec:shadow}

As we have discussed, numerical trajectories and exact trajectories started at the same initial condition typically rapidly diverge.
One idea that has been proposed is that for every numerical trajectory there is an exact trajectory starting at a different initial condition that stays close to the numerical trajectory over long time intervals. 
  This idea is known as shadowing.
An early  version of shadowing is described by Bowen \cite{bowen} for Axiom A systems, though see \cite{anosovshadow}, \cite[p. 381]{bowenearly}, and  \cite[p.\ 38]{sinaishadow} for earlier descriptions.  A general result  in this area is that  if a system satisfies a uniform hyperbolicity condition, then shadowing is possible over infinite time intervals \cite{pilyugin}.   This fact was used in \cite{Reich} to study the long-time averages over trajectories computed with a symplectic method under the assumption that the Poincar\'e section of the flow is uniformly hyperbolic.

However, many systems that arise in applications are not uniformly hyperbolic \cite[Appendix B]{hunt}.  To the best of our knowledge, the only example of a physically realizable Hamiltonian system that is uniformly hyperbolic on one of its energy levels is due to Hunt and Mackey \cite{hunt}, and this system is uncharacteristic of systems that arise in molecular dynamics.  (Many billiard systems have been shown to be ergodic and even mixing, but fail to be uniformly hyperbolic because the vector field is discontinuous at bounces \cite[Appendix B]{hunt}.)
  For more realistic systems shadowing has been numerically demonstrated over long but finite time intervals \cite{saueryorke}, \cite{hayesprl}.    It remains to be seen whether shadowing over the long trajectories computed in molecular dynamics is possible.

Let us specify formally what shadowing would consist of in our situation.  Fixing a time interval $[0,T]$ we say that $Y$, an actual trajectory of the system, $\epsilon$-shadows the numerical trajectory $X_{\Delta t}$ if
\[
\| Y - X_{\Delta t} \|_\infty 
\leq \epsilon.
\]
Assuming that it is possible to shadow every numerical trajectory in
this way, let $\mathcal{S}_{\Delta t}$ be the map that takes the initial condition of the numerical trajectory $X_{\Delta t}$ to the initial condition of the exact trajectory $Y$ that shadows $X_{\Delta t}$.  This gives us our first version of shadowing.

{\bf Shadowing Version 1}. 
\emph{ There exists $\mathcal{S}_{\Delta t} \colon \mathbb{R}^m \rightarrow \mathbb{R}^m$ such that if $Y_0=\mathcal{S}_{\Delta t} X_0$ then }
\[
\| Y - X_{\Delta t} \|_\infty  \leq \epsilon.
\]

This version of shadowing is not sufficient for our purposes: it is not strong enough to  show something like  (\ref{eqn:goal}).   The difficulty is that even though each numerical trajectory is close to some exact trajectory, it could be that the distribution of $Y$ is completely different from that of $X$.  This can result even if the distribution of $X_0$ and $Y_0$ are close because of the rapid divergence of trajectories of the system.  This would in turn imply, since $X_{\Delta t}$ is close to $Y$ for every $X_0$, that $X_{\Delta t}$ is not close to $X$ in distribution.

In order for the statistical properties of $Y$ to be similar to that of $X$, it is necessary that the map $\mathcal{S}_{\Delta t}$ preserves the measure on the initial condition $X_0$.  We are lead to the following concept of shadowing for initial value problems with a probability measure on the initial conditions:\\

{\bf Shadowing Version 2}. \emph{Given that $X_0$ has distribution $\nu$ on $\mathbb{R}^m$ there is a map $\mathcal{S}_{\Delta t}$ on $\mathbb{R}^m$ such that $\mathcal{S}_{\Delta t} X_0$ also has distribution $\nu$ and}
 \[
 \| Y - X_{\Delta t} \|_{\infty}  \leq \epsilon.
\]

In practice we want to allow the possibility that shadowing is not possible for some 
 initial conditions.  So to weaken Shadowing Version 2 slightly we conjecture:

{\bf Shadowing Version 3.}  \emph{Given that $X_0$ has distribution $\nu$ on $\mathbb{R}^m$ there is a map $\mathcal{S}_{\Delta t}$ on $\mathbb{R}^m$ such that $\mathcal{S}_{\Delta t} X_0$ also has distribution $\nu$ and }

\[
\mathbb{P}_{\nu} ( \| Y - X_{\Delta t} \|_{\infty}  > \epsilon  ) < \epsilon.
\]

With our application in mind there is an important way we can further weaken Shadowing Version 3.  
As in the previous subsection, 
let $\Pi \colon \mathbb{R}^m \rightarrow \mathbb{R}^k$ be a map that takes configurations of the whole system and extracts lower dimensional information.
 Then we could conjecture:

 {\bf Shadowing Version 4.} \emph{Given that $X_0$ has distribution $\nu$ on $\mathbb{R}^m$ there 
 is a map $\mathcal{S}_{\Delta t}$ on $\mathbb{R}^m$ 
 such that $\mathcal{S}_{\Delta t} X_0$ also has distribution $\nu$ and
 }
\begin{equation} \label{eqn:sv4}
\mathbb{P}_{\nu} ( \| \Pi(Y) - \Pi(X_{\Delta t}) \|_{\infty}  > \epsilon  ) < \epsilon.
\end{equation}

This final form of shadowing is what we call \emph{weak shadowing}.  

In the context of molecular dynamics, establishing weak shadowing for some physically interesting $\Pi$ would be relevant  if for all sufficiently small $\Delta t$ there were constants $A,B,C,p>0$ such that  (\ref{eqn:sv4}) held with $\epsilon = A \Delta t^p$ for all $T = B \exp(-C/\Delta t)$.
Note that all versions of shadowing presented here hold immediately for small $\epsilon$ if $T$ is fixed and $\Delta t$ is allowed to be arbitrarily small.  However, this limit is not of interest in molecular dynamics.

In order for weak shadowing to have the same explanatory power as Approximation in Distribution, we need to show that it also allows us to bound the difference between $\mathbb{E} G(X)$ and $\mathbb{E} G(X_{\Delta t})$ for reasonable functions $G$, as in (\ref{eqn:goal}).  To see that it does, again define $\| \cdot \|_{BL}$ as in the previous subsection and note that for any $G \colon (C[0,T])^k \rightarrow \mathbb{R}$, $G$ is bounded by $\| G\|_{BL}$ and $G$ has $\| G\|_{BL}$ as a Lipschitz constant.
In that case:
\begin{eqnarray*}
| \mathbb{E} G(\Pi(X)) - \mathbb{E} G( \Pi(X_{\Delta t})) |&  = &
| \mathbb{E} G(\Pi(Y)) - \mathbb{E} G( \Pi(X_{\Delta t})) | \\
& \leq & \|G\|_{BL} \epsilon + 2 \|G\|_{GL} \mathbb{P}( \| \Pi(Y) - \Pi(X_{\Delta t}) \| > \epsilon ) \\
& < & 3 \| G\|_{BL} \epsilon,
\end{eqnarray*}
where we have used the fact that $X$ and $Y$ have the same distribution in the first equality.
This result is precisely (\ref{eqn:goal}) with a different constant.

\section{Proof of Main Result}    \label{sec:mainthm}

Our main result, Theorem~\ref{thm:easythm} in Section~\ref{sec:intro},  asserts that Approximation in Distribution (Equation (\ref{eqn:approxdist})) and Weak Shadowing (Equation (\ref{eqn:sv4})) are equivalent.  In this section we first state a more general result, Theorem~\ref{thm:mainthm}, and then show how Theorem~\ref{thm:easythm} follows from it.  Then we present the proof of Theorem~\ref{thm:mainthm} which uses our main technical result, Theorem~\ref{thm:maintechnical}.  

We say that a measure $\nu$ on a space $S$ is \emph{non-atomic} if $\nu(\{x\}) =0$ for every point $x \in S$.

\begin{theorem} \label{thm:mainthm}
Let  $\Phi$, $\hat{\Phi} \colon \mathbb{R}^m \rightarrow (C[0,T])^k$ be measurable maps and let $\nu$ be a non-atomic measure on $\mathbb{R}^m$.  Let $X_0$ be a random vector in $\mathbb{R}^k$ with distribution $\nu$. Then the following are equivalent for all $\epsilon>0$: \\
(A) {\bf Approximation in distribution:}
\[
\rho(\Phi(X_0),\hat{\Phi}(X_0) )  < \epsilon.
\]
(B) {\bf Weak Shadowing:}  There is a map $\mathcal{S} \colon \mathbb{R}^m \rightarrow \mathbb{R}^m$ 
such that $\mathcal{S}X_0$ also has distribution $\nu$ and 
\[
\mathbb{P}(  \|  \Phi(\mathcal{S} X_0) - \hat{\Phi}(X_0)   \|_{\infty} > \epsilon) < \epsilon.
\]
\end{theorem}
{\em Proof}.
To show that (A) implies (B) it is only necessary to apply Theorem~\ref{thm:maintechnical} with $(S,d)= (T,e)= (\mathbb{R}^m,\| \cdot \|)$, $X=Y=X_0$, $(\bar{S},\bar{d}) = ((C[0,T])^k,\| \cdot \|_\infty)$, $\Pi_1= \Phi_1$, $\Pi_2=\Phi_2$.
   The theorem then gives a map $\psi \colon \mathbb{R}^m \rightarrow \mathbb{R}^m$ such that $\psi X_0$ has the same distribution as $X_0$ and 
\[
\mathbb{P}( \| \Phi(\psi(X_0))) - \hat{\Phi}(X_0) \|_\infty > \epsilon )< \epsilon.
\]

To show that (B) implies (A), it is only necessary to see that (B) implies
\[
\rho( \Phi(\mathcal{S} X_0) , \hat{\Phi}(X_0) ) <\epsilon
\]
 Since $\mathcal{S} X_0$ is equal in distribution to $X_0$ and equality in distribution is preserved under measurable maps, $\Phi(\mathcal{S} X_0)$ is equal to $\Pi(X_0)$ in distribution.
Since the Prokhorov metric $\rho$ only depends on distributions, we have that $\rho(\Phi(X_0), \hat{\Phi}(X_0))<\epsilon$ as required.
$\endproof$

{\em Proof of Theorem~\ref{thm:easythm}.}
Let $\Phi$ be the map that takes initial condition $X_0$ to  $\Pi(X)$.
Let $\hat{\Phi}$ be the map that takes $X_0$ to $\Pi(X_{\Delta})$.  Let $\nu$ be the distribution of $X_0$.  Then the theorem follows.
$\endproof$

{\bf Remark:} Though we have motivated our result in terms of $\Phi$ being  the exact flow of a system of differential equations and $\hat{\Phi}$ being the trajectory generated by a numerical method, the result can be much more broadly applied.  In particular, $\Phi$ and $\hat{\Phi}$ can be the flow maps of any two dynamical systems on the same state space.  

It only remains to prove Theorem~\ref{thm:maintechnical}, which
 is the heart of Theorem~\ref{thm:mainthm} above.
Theorem~\ref{thm:maintechnical} is  similar to the Strassen-Dudley Theorem \cite{dudley,strassen} which we state here for reference:
\begin{theorem} \label{thm:strass} (\cite[p.\ 73]{billingsley})  Let $(S,d)$ be a separable metric space.
 If $X$ and $Y$ are random elements of $S$ with 
  $\rho(X,Y) < \beta$, then there are random elements $\bar{X}$ and $\bar{Y}$ of $S$ defined on a common probability space
    such that  $\bar{X}$ has the same distribution as $X$, $\bar{Y}$ has the same distribution as $Y$ and 
\[
\mathbb{P} ( d(\bar{X},\bar{Y}) > \beta ) < \beta. 
\]
\hfill $\square$
\end{theorem}

In contrast, our theorem is as follows.  Recall that the distribution of $X$, a random element of $S$, is non-atomic when $\mathbb{P}( X = x) =0$ for all $x \in S$.
\begin{theorem} \label{thm:maintechnical}
Let $X$ be a random element of the separable complete metric space $(S,d)$.  Let $Y$ be a random element of the separable complete metric space $(T,e)$.  Suppose that the distributions of both $X$ and $Y$ are non-atomic.
 Let $(\bar{S},\bar{d})$  be another separable complete metric space, and let $\Pi_1 \colon S \rightarrow \bar{S}$ and $\Pi_2 \colon T \rightarrow \bar{S}$ be measurable maps.    Let $\rho$ denote the Prokhorov metric on random elements of $(\bar{S},\bar{d})$.  If $\rho(\Pi_1(X),\Pi_2(Y)) < \beta$
then there is a measurable map $\psi: S \rightarrow T$ such that $\bar{Y} = \psi X$ has the same distribution as $Y$ and 
\begin{equation} \label{eq:techthmcond}
\mathbb{P} ( \bar{d}(\Pi_1(X),\Pi_2(\bar{Y})) > \beta ) < \beta.
\end{equation}
\end{theorem}

Taking $(\bar{S},\bar{d})= (S,d)=(T,e)$ and  $\Pi_1$ and $\Pi_2$ to be the identity gives the following simple corollary that is easier to compare with Theorem~\ref{thm:strass}.

\begin{corollary} \label{cor:conseq}  Let $(S,d)$ be a separable complete metric space.  Let $X$ and $Y$ be random elements of $S$ with non-atomic distributions.
If $\rho(X,Y) < \beta$, then there is measurable map $\psi$ from $(S,d)$ to itself such that 
$\bar{Y} = \psi X$ has the same distribution as $Y$ and
\[
\mathbb{P}(d(X,\bar{Y})>\beta) < \beta.
\]
\end{corollary}

This result is stronger than Theorem~\ref{thm:strass} in that instead of being forced to have a new probability space and define two new random variables $\bar{X}$ and $\bar{Y}$, we are able to leave $X$ as it is and define $\bar{Y}$ to be a random variable on the same space that $X$ is defined on.  This extra strength is necessary for us in order to make the connection with shadowing in Theorem~\ref{thm:mainthm}.  The extra cost is that we assume that  the metric spaces in which $X$ and $Y$ live are complete and, more importantly, that the probability distributions of $X$ and $Y$ are non-atomic.

To see that the assumption of non-atomicity is essential in Corollary~\ref{cor:conseq}, and hence in Theorem~\ref{thm:maintechnical}, consider the following example.  For any  $\epsilon \in (0,1/2)$ let $\Omega = \left\{ 0, 1 \right\}$ and let $\mathcal{F}$ be all the subsets of $\Omega$.  Let the probability measure $\mathbb{P}$ be defined by $\mathbb{P}(0)=1/2-\epsilon$, $\mathbb{P}(1) = 1/2+\epsilon$.  Then $(\Omega,\mathcal{F},\mathbb{P})$ is a probability space.  Define random variables $X$ and $Y$ by $X(0)=0, X(1)=1, Y(0)=1, Y(1)=0$.  It is straightforward to check that $\rho(X,Y)=\epsilon$.  However, $Y$ is the \emph{only} random variable on 
 $(\Omega, \mathcal{F}, \mathbb{P})$ that has the same distribution as $Y$, and 
 $\mathbb{P} ( |X-Y| > 1/2 ) = 2 \epsilon$.  So $\mathbb{P}( |X-Y| > \epsilon ) > \epsilon$, and the result cannot hold.

  Now we turn to the proof of Theorem~\ref{thm:maintechnical}.
  Many of the ideas used in the proof of the Strassen-Dudley Theorem reappear here, including the use of the Marriage Lemma.  Before the proof we need several lemmas that allow us to construct measure isomorphisms between various spaces and also to divide spaces into many pieces of equal measure.

\begin{lemma} \label{lem:fromroyden} \cite[p. 327]{Royden}
Let $(S,d)$ be a complete separable metric space with Borel $\sigma$-algebra $\mathcal{B}(S)$.  Let $\mu$ be a non-atomic probability measure on $(S, \mathcal{B}(S))$.   Let $([0,1],$ $\mathcal{B}([0,1]),$ $\lambda)$ be the unit interval with Lebesgue measure defined on the Borel sets.
Then there is a subset $S_0$ of $S$ with $\mu(S_0)=\mu(S)$ and subset $L_0$ of $[0,1]$ with $\lambda(L_0) = \lambda([0,1])$ such that there is a measurable invertible $\psi : S_0 \rightarrow L_0$ such that $\mu(\psi^{-1} A) = \lambda (A)$ for all $A \in \mathcal{B}([0,1]) \cap L_0$. $\endproof$
\end{lemma}

\begin{lemma} \label{cor:fleshout}
Let $(S,d)$ and $(T,e)$ be two complete separable metric spaces with Borel $\sigma$-algebras $\mathcal{B}(S)$ and $\mathcal{B}(T)$ respectively.  Let $\mu_S$ and $\mu_T$ be non-atomic probability measures on $(S,\mathcal{B}(S))$ and $(T,\mathcal{B}(T))$ respectively.   Then there is a measurable map $\psi \colon S \rightarrow T$ such that
\[
\mu_S(\psi^{-1} A) = \mu_T( A)
\]
for all $A \in \mathcal{B}(T)$.
\end{lemma}
{\em Proof}.  Use Lemma \ref{lem:fromroyden} to construct subsets of full measure $S_0$ and $T_0$ in $S$ and $T$ with measure preserving invertible maps $\psi_1: S_0\rightarrow L_1$ and $\psi_2: T_0 \rightarrow L_2$.  Let $\tilde{S}_0 = \psi_1^{-1} L_1 \cap \psi_1^{-1} L_2$  and let $\psi$ equal $\psi_2^{-1} \psi_1$ restricted to this set.  Let $\tilde{T}_0$ be the image of $S_0$ under $\psi$.  Now extend $\psi$ in an arbitrary measurable way to all of $S$.
$\endproof$

\begin{lemma} \label{cor:dice}
Let $(S,d)$ be a complete separable metric space with Borel $\sigma$-algebra $\mathcal{B}(S)$.  Let $\mu$ be a non-atomic probability measure on $(S,\mathcal{B})$.   Given any subset $T \in \mathcal{B}(S)$ with  $\mu(T) = m \leq 1$ and any finite set of real numbers $m_i$, $i=1,\ldots,n$ with $\sum_i m_i = m$ there is a partition of $T$
\[
T = \cup_{i=1}^n T_i
\]
with $\mu(T_i) = m_i$, $i= 1,\ldots , n$.
\end{lemma}

{\em Proof}.
We begin by proving the result for the $(S,d)$ being the unit interval $[0,1]$ and $\mu$ Lebesgue measure.  Then we use the previous lemma to establish the general case.

Let $T$ be a Borel subset of $[0,1]$ with Lebesgue measure $m$.  Consider the function $\sigma: [0,1] \rightarrow [0,m]$ defined by 
\[
\sigma(x) = \lambda ( T \cap [0,x]).
\]
Since $\lambda$ is non-atomic $\sigma$ is continuous.
Let $\bar{m}_0=0$ and $\bar{m}_i = \sum_{j=1}^i m_j$.
Let $T_i = \sigma^{-1} ( [\bar{m}_{i-1}, \bar{m}_{i}))$ for $1 \leq i \leq n-1$ and $T_n = \sigma^{-1} ([ \bar{m}_{n-1} \bar{m}_n])$.  It is straightforward to show that the $T_i$ satisfy the required conditions.

For the case of general $(S,d)$, let $\psi$, $S_0$, $L_0$ be as given by Lemma~\ref{lem:fromroyden}.  The subset of $[0,1]$ given by $\psi(T \cap S_0)$ has measure $m$.  By the previous case this can be partitioned into $n$ subsets $R_i$ with $\lambda(R_i) = m_i$ that are all subsets of $L_0$.  Let $T_i = \psi^{-1}(R_i)$ for $1 \leq i \leq n-1$ and $T_n = \psi^{-1}(R_n) \cup (T \setminus S_0)$.  It is straightforward to show that $T_i$ have the required properties.
$\endproof$

\begin{lemma} \label{lem:onepart}
Let $(S,d)$ be a complete separable metric space with non-atomic  probability measure $\mu$ on it.  Let $(\bar{S},\bar{d})$ be another complete separable metric space and $\Pi \colon S \rightarrow \bar{S}$ a measurable map.   For any $\epsilon>0$ there is a $\delta < \epsilon$  and a finite partition of $S$
\[
S= S^* \cup ( \cup_{k=1}^n S_k),
\]
such that \\
(i) $\mu(S_k) = \delta$ for all $k$,\\
(ii) $\mathrm{diam}(\Pi(S_k)) < \epsilon$ for all $k$, \\
(iii) $\mu(S^*)< \epsilon$. \\
Moreover, this is possible for all sufficiently small $\delta$.
\end{lemma}

{\em Proof}.
Let $x_i$, $i \geq 1$ be a dense sequence of points in $\bar{S}$.  Let $B_i \subset S$ be the inverse image under $\Pi$ of the open ball of radius $\epsilon/2$ about $x_i$ in $\bar{S}$.  For $i\geq 1$ let
\[
\bar{B}_i = B_i  \setminus \cup_{j=1}^{i-1} B_j.
\]
Then the $\bar{B}_i$ are disjoint, $\mathrm{diam} (\Pi (\bar{B}_i)) < \epsilon$ and have union $S$.  Choose $m$ such that 
\[
\sum_{i=1}^m \mu(\bar{B}_i) > 1-\epsilon/2. 
\]
Now choose a  sufficiently small $\delta$ so that $m \delta < \epsilon/2$.  Using Corollary~\ref{cor:dice} divide each $\bar{B}_i$ into $k_i$ sets $\bar{B}_{i,j}$ with $\mu(\bar{B}_{i,j}) = \delta$ and one additional set $\bar{B}_i^*$ with $\mu(\bar{B}_i^*)<\delta$.  Now there are finitely many sets  $\bar{B}_{i,j}$ all  having measure $\delta$.  Call these sets $S_k$.  Condition (i) is then satisfied. Each has diameter less than $\epsilon$, since each is a subset of some $\bar{B}_i$.   So condition (ii) is satisfied.  If we let  $S^*= S \setminus \cup_{k=1}^n S_k$ 
we have 
\begin{eqnarray*}
\mu( S^* ) & =&  \mu \left( S \setminus \sum_{i=1}^m \mu(\bar{B}_i)\right) + \sum_{i=1}^m \mu(\bar{B}_i^*) \\
   & \leq &  \epsilon/2 + m \delta =  \epsilon,
\end{eqnarray*}
showing that condition (iii) is satisfied.
$\endproof$

\emph{Proof of Theorem~\ref{thm:maintechnical}.}  Let $\alpha = \rho(\Pi_1(X),\Pi_2(Y)) < \beta$.
Consider any $\epsilon>0$, which we will fix later to obtain the required result.

Use Lemma~\ref{lem:onepart} to construct  finite partitions of $S$ and $T$,
\[
S= S^* \cup ( \cup_{i=1}^n S_i), \ \ \ \ T = T^* \cup( \cup_{i=1}^n T_i), 
\]
such that
\[
 \mathbb{P}(X \in S^*)< \epsilon, \ \ \ \ \ \  \mathbb{P}( Y \in T^*)< \epsilon,
 \]
 and  for all $i$ 
\[
\mbox{diam}(\Pi_1(S_i)) < \epsilon,  \ \ \ \ \ \mbox{diam}(\Pi_2(T_i)) < \epsilon,  
\]
\[
 \mathbb{P}(X \in S_i) = \mathbb{P}(Y \in T_i) = \delta < \epsilon.
 \]
We can use the same $\delta$ for both $S$ and $T$ since Lemma~\ref{lem:onepart} shows that for each $\epsilon$ the construction is possible for all sufficiently small $\delta$.

We will construct a 1-1 mapping $\phi$  from the set $\{1,\ldots,n\}$ to itself such that for most $i$ we have
$\bar{d}(\Pi_1(S_i), \Pi_2(T_{\phi(i)}))< \alpha + \epsilon$.   In other words, for most $i$ there will be a point in $\Pi_1(S_i)$ and a point in $\Pi_2(T_{\phi(i)})$ that are within distance $\alpha+\epsilon$ of each other.  
Based on $\phi$,  we will then use Lemma~\ref{lem:fromroyden} to construct a map $\psi$ on $S$ that takes $S_i$ to $T_{\phi(i)}$ and such that $\psi X$ has the same distribution as $Y$.  We will construct the map $\phi$ with the help of the Marriage Lemma of K\"onig and Hall \cite{dudley}.

\begin{lemma} (See \cite[p.\ 406]{dudley}.)
Let $K$ denote a relation on $\{1,\ldots,n\}$ such that for all subsets $A$ of $\{1,\ldots,n\}$
\begin{equation}  \label{eqn:card}
| \{ j \in A  \colon {}_iK_j \mbox{ for some } i \in A \} | \geq | A|
\end{equation}
where $| \cdot |$ denotes cardinality.  Then there is a 1-1 mapping $\phi$ of $\{1,\ldots,n\}$ to itself such that ${}_iK_{\phi(i)}$ for all $i$. $\endproof$
\end{lemma}

Ideally we would define the relation $K$ on $\{ 1, \ldots, n\}$
by saying that ${}_iK_j$ if $\bar{d}(\Pi_1(S_i),\Pi_2(T_j))$ $< \alpha+ \epsilon$, and then using the Marriage Lemma to construct a mapping $\phi$ such that $\bar{d}(\Pi_1(S_i),\Pi_2(T_{\phi(i)}))< \alpha +\epsilon$ for all $i$.  
However, in general \eqref{eqn:card} does not hold for this definition of $K$, and such a $\phi$ does not exist.

Instead, we construct a map $\phi$ such that $\bar{d}(\Pi_1(S_i),\Pi_2(T_{\phi(i)}))< \alpha+\epsilon$ only for most $i$, as follows.
We append to $\{1, \ldots, n\}$ $k$ extra indices $n+1, \ldots, n+k$.  Now for $i,j \in \{ 1, \ldots, n+k \}$ we say that ${}_iK_j$ if either
\begin{enumerate}
\item   $\bar{d}(\Pi_1(S_i),\Pi_2(T_j)) < \alpha +\epsilon$,
\item  $i > n$, or
\item  $j > n$.
\end{enumerate}

Now let $A$ be a subset of $\{1, \ldots, n+k \}$.  Either $A$ contains at least one of $n+1, \ldots, n+k$ or it doesn't.    In the former case \eqref{eqn:card} holds immediately.  In the latter case,  let $S_A = \cup_{i \in A} S_i$.   We have that 
\[
\mathbb{P}(X \in S_A) = | A | \delta.
\]  
Let 
\[
B = \{ z \in S \colon \bar{d}(\Pi_1(S_A),\Pi_2(z)) < \alpha + \epsilon \},
\]
so that 
\begin{eqnarray*}
\mathbb{P}(Y \in B)  & = & \mathbb{P} ( \bar{d}(\Pi_1(S_A),\Pi_2(Y)) < \alpha +\epsilon )\\
& \geq & \mathbb{P} ( \bar{d}(\Pi_1(S_A),\Pi_2(Y)) \leq \alpha ) \\
& = & \mathbb{P} (\Pi_2(Y) \in \Pi_1(S_A)^\alpha )
\end{eqnarray*}
Then since $\rho(\Pi_1(X),\Pi_2(Y))= \alpha$,
\[
\mathbb{P} ( \Pi_1(X) \in \Pi_1(S_A) ) \leq \mathbb{P} ( \Pi_2(Y) \in \Pi_1(S_A)^\alpha ) + \alpha.
\] 
This fact yields
\[
 \mathbb{P}(Y \in B)  \geq \mathbb{P}( \Pi_1(X) \in \Pi_1(S_A) ) - \alpha 
 \geq \mathbb{P}(X \in S_A) - \alpha=   |A| \delta -\alpha.
\]
So the number of sets $T_j$ that have some portion in $B$ is at least
\[
(|A| \delta -\alpha - \mathbb{P}(Y \in T^* ))/ \delta 
\geq |A| -(\alpha + \epsilon)/\delta. 
\]
For all these $j \leq n$  there is some $i \in A$ such that ${}_iK_j$.

When we include all the $j > n$,
the total number of $j$ such that ${}_iK_j$ for some $i \in A$ is then at least
$|A|  -(\alpha+\epsilon)/\delta + k$.
So if we let $k =\lceil (\alpha+\epsilon)/\delta \rceil$, and then relation $K$ on $\{ 1, \ldots, n+k \}$ satisfies the conditions of the Marriage Lemma.

This gives us that there is a 1-1 map $\bar{\phi}$  between $\{ 1, \ldots, n+k \}$ and itself such that ${}_i K_{\bar{\phi}(i)}$ for all $i$.  
We want to get a map 1-1 map $\phi$ on $\{ 1, \ldots, n \}$  so that ${}_iK_{\phi(i)}$ for \emph{most} $i$.  
Consider what $\bar{\phi}$ does to the set $\{ 1, \ldots, n \}$.  At least $n-k$ elements get mapped back to $\{ 1, \ldots, n \}$.  Let $\phi(i)= \bar{\phi}(i)$ for these elements.  For all the others, just let $\phi(i)$ be extended to be 1-1 on $\{1,\ldots, n\}$.

Now we have an invertible map $\phi$ on $\{1, \ldots, n\}$ such that for $n-k$ of the elements $i$
\[
\bar{d}( \Pi_1(S_i), \Pi_2(T_{\phi(i)})) < \alpha+ \epsilon.
\] 
Using Lemma~\ref{cor:fleshout}, for each $i$ there is a map $\psi_i : S_i \rightarrow T_{\phi(i)}$ such that for any measurable subset $C$ of $S_i$
\[
\mathbb{P}( X \in C ) = \mathbb{P} (Y \in \psi_i(C)).
\]
Again use Lemma~\ref{cor:fleshout} to construct a map $\psi_* \colon S^* \rightarrow T^*$ such that $\mathbb{P}( X \in C ) = \mathbb{P} (Y \in \psi_*(C))$ for all measurable $C \subset S^*$.
Now let $\psi \colon S \rightarrow T$ be defined by requiring that $\psi$ restricted to $S_i$ is $\psi_i$ and that $\psi$ restricted to $S^*$ is $\psi_*$.  Then for any measurable subset  $C$ of $S$
\[
\mathbb{P} (X \in C) =  \mathbb{P} ( Y \in \psi C).
\]
This means that if we let $\bar{Y}=\psi X$, for any measurable subset $D$ of $S$
\[
\mathbb{P}( \bar{Y} \in D ) = \mathbb{P}( \psi X \in D) = \mathbb{P}(X \in \psi^{-1} D) = 
  \mathbb{P}( Y \in  \psi \psi^{-1} D ) = \mathbb{P} (Y \in D)
\]
as required.

It remains to show that $\bar{Y}$ satisfies equation (\ref{eq:techthmcond}).
Now for $n-k$ indices $i$
\[
\bar{d} ( \Pi_1(S_i) , \Pi_2(T_{\phi(i)}) ) < \alpha + \epsilon.
\]
If $X \in S_{i}$ then $\bar{Y} \in T_{\phi(i)}$ and 
\[
\bar{d}( \Pi_1(X) , \Pi_2(\bar{Y}) )< \alpha + \epsilon + 2\epsilon,
\]
since $\Pi_1(S_i)$ and $\Pi_2( T_{\phi(i)})$ have diameters smaller than $\epsilon$.
So with probability at least  $\delta (n-k)$ we have that 
$ \bar{d}( \Pi_1(X) , \Pi_2(\bar{Y}) )< \alpha + 3 \epsilon$.
Since $1 = n \delta + \mathbb{P}(X \in S^*)$ and  $\mathbb{P}(X \in S^*) < \epsilon$,
\[
\mathbb{P}( \bar{d}( \Pi_1 (X), \Pi_2( \bar{Y} ) ) \geq \alpha +3\epsilon ) < k \delta + \epsilon =  \delta \lceil (\alpha + \epsilon)/ \delta \rceil   + \delta \leq \alpha + \epsilon +2 \delta < \alpha + 3\epsilon.
\]
Choosing $\epsilon$ so that $\alpha + 3 \epsilon <\beta$ then gives our result.
\hfill$\square$

\section{Discussion} \label{sec:discussion}

Suppose we are considering a particular molecular dynamics simulation over a long time interval $[0,T]$ started from random initial conditions.
We wish to determine what statistical features of its 
trajectories are computed reliably.   A simple baseline conjecture would be that all statistical features of the trajectories are computed accurately.  As we have detailed above, a quantitative version of this conjecture would be to say that 
$\rho(X,X_{\Delta t}) =\epsilon$ for some small $\epsilon > 0$.   Then Theorem~\ref{thm:easythm} states that with probability greater than $1-\epsilon$, 
 the numerical trajectory is shadowed by an exact trajectory to within error $\epsilon$.
In this case, it should be possible to detect shadow trajectories numerically using the techniques of \cite{hayessisc}, even though the shadows computed will not necessarily have the correct measure on their initial conditions.  Find such shadow trajectories would be partial confirmation that $\rho(X,X_{\Delta t})$ is small.

The other possibility is that $\rho(X,X_{\Delta t})$ is not small.  
Suppose instead that $\rho(X,X_{\Delta t})>1/2$.  Now (\ref{eqn:rhobeta}) implies that $\beta(X,X_{\Delta t}) > 1/6$.  This means that 
\begin{equation}  \label{eqn:problem}
\| \mathbb{E}G(X) - \mathbb{E}G(X_{\Delta t}) \| > 1/6 
\end{equation}
for some $G \colon (C[0,T])^m \rightarrow \mathbb{R}$ with $\|G\|_{BL}$.
Hence one way to confirm that $\rho(X,X_{\Delta t})$ is large is to find such a $G$ such that we empirically observe (\ref{eqn:problem}).

 In principle this approach is reasonable, but in practice the space of all functions $G \colon (C[0,T])^m \rightarrow \mathbb{R}$ with $\| G \|_{BL}=1$ is huge for a realistic molecular dynamics simulation.   In practice it may be that only for very unusual functions $G$ do $\mathbb{E}G(X)$ and $\mathbb{E}G(X_{\Delta t})$ disagree significantly.  
One practical way to approach this is to start with very low-dimensional systems.  The lowest dimensional system that is a reasonable model of molecular dynamics consists of two particles on a two-dimensional periodic domain \cite{szasz}.  For example, one could study the system we considered in Section~\ref{sec:problem} but with only two particles.   Using the software available to compute exact shadow trajectories of numerical trajectories \cite{hayessisc} would show where shadowing is not possible and could suggest what functions $G$ are likely candidates.
 
We have been considering the case where $\Pi$ is the identity, for which it may be that $\rho(\Pi(X),\Pi(X_{\Delta t}))$ is large.  The other direction to study these systems is to choose a $\Pi$ that is a very low dimensional function of the state of the system and then see if it is possible to numerically perform weak shadowing with this choice.  Currently there have not been algorithms developed to do this,  but it is a direction for future work.

Finally, besides these numerical/experimental approaches there are more analytical approaches.  These would involve studying one of the model systems for molecular dynamics available that are analytically tractable.  Both the systems studied in \cite{KSTT} and \cite{Murmann} consist of single particle coupled to a bath of a very large or infinite number of smaller particles.  In both cases it is shown that the distribution of the path of the large particle converges to a stochastic process that can be simply described through low-dimensional stochastic differential equations.  These situations are clean enough that similar results for the numerical discretization of the Hamiltonian equation may be possible, thus establishing approximation in distribution for the path of the distinguished particle.  Our result then allows us to conclude that weak shadowing is possible.

Our work allows the possibility of studying the reliability molecular dynamics in a variety of contexts from one of two directions: either the statistical one through the computing of histograms, or the dynamical one through the computing of shadowing trajectories.
 
{\bf Acknowledgements.}  The author would like to thank Robert Skeel, Wayne Hayes, and Nilima Nigam for comments on earlier versions of this manuscript.

\bibliography{../../paulsbib}{}
\bibliographystyle{abbrv}

\end{document}